\documentclass[times]{elsarticle}

\usepackage[dvipsnames]{xcolor}
\usepackage{amsmath}
\usepackage{amsfonts}
\usepackage{amssymb}
\usepackage{lineno}
\usepackage{enumerate}
\usepackage{times}
\usepackage{subcaption}
\usepackage{graphicx,psfrag}
\usepackage{pgfplotstable}
\usepackage[skip=0pt]{caption}

\newcommand\solidrule[1][0.25cm]{\rule[0.5ex]{#1}{1pt}}
\newcommand\dashedrule{\mbox{%
  \solidrule[2mm]\hspace{2mm}\solidrule[2mm]}}

\newcommand{\dotrule}[1]{%
	\parbox{#1}{\dotfill}} 

\makeatletter
\newcommand \Dotfill {\leavevmode \cleaders \hb@xt@ .22em{\hss .\hss }\hfill \kern \z@}
\makeatother

\begin{document}

\title{Behaviour of the Serre Equations in the Presence of Steep Gradients Revisited}

\author[ANU]{J.P.A.~Pitt\corref{cor1}}
\ead{jordan.pitt@anu.edu.au }
\author[ANU]{C.~Zoppou}
\ead{christopher.zoppou@anu.edu.au}
\author[ANU]{S.G.~Roberts}
\ead{stephen.roberts@anu.edu.au}

\cortext[cor1]{Corresponding author}
\address[ANU]{Mathematical Sciences Institute, Australian National University, Canberra, ACT 0200, Australia}
 \begin{abstract}
 We use numerical methods to study the behaviour of the Serre equations in the presence of steep gradients because there are no known analytical solutions for these problems. In keeping with the literature we study a class of initial condition problems that are a smooth approximation to the initial conditions of the dam-break problem. This class of initial condition problems allow us to observe the behaviour of the Serre equations with varying steepness of the initial conditions. The numerical solutions of the Serre equations are justified by demonstrating that as the resolution increases they converge to a solution with little error in conservation of mass, momentum and energy independent of the numerical method. We observe four different structures of the converged numerical solutions depending on the steepness of the initial conditions. Two of these structures were observed in the literature, with the other two not being commonly found in the literature. The numerical solutions are then used to assess how well the analytical solution of the shallow water wave equations captures the mean behaviour of the solution of the Serre equations for the dam-break problem. Lastly the numerical solutions are used to evaluate the usefulness of asymptotic results in the literature to approximate the depth and location of the front of an undular bore.
 \end{abstract}	
 
  
 \maketitle
\section{Introduction} \label{intro} 
The behaviour of flows containing steep gradients are important to a range of problems in shallow water such as the propagation of a bore, the dam-break problem and shoaling waves on a beach.

The Serre equations are used as a compromise between the non-dispersive shallow water wave equations and the incompressible inviscid Euler equations for modelling dispersive waves of the free surface in the presence of steep gradients, which are present for the Euler equations \cite{Mitsotakis-etal-2017} but not for the shallow water wave equations. The Serre equations like the shallow water wave equations produce methods \cite{Hank-etal-2010-2034,Zoppou-etal-2016,Zoppou-etal-2017} that are computationally easier and quicker to solve than the best methods for the Euler equations. The Serre equations are considered the most appropriate approximation to the Euler equations for modelling dispersive waves up to the shore line \cite{Bonneton-etal-2011-589,Bonneton-etal-2011-1479}. Therefore, understanding the behaviour of the Serre equations in the presence of steep gradients offers some insight into the behaviour of steep gradients for fluids more generally.

There are no known analytical solutions to problems containing steep gradients for the Serre equations. To infer the structure of solutions to problems containing steep gradients we have to resort to investigating numerical solutions of the Serre equations for these problems. 

There are few examples in the literature which depict the behaviour of numerical solutions to the Serre equations in the presence of steep gradients \cite{El-etal-2006,Hank-etal-2010-2034,Mitsotakis-etal-2014,Mitsotakis-etal-2017,Zoppou-etal-2016,Zoppou-etal-2017}. These papers all present problems with discontinuous initial conditions \cite{Hank-etal-2010-2034,Zoppou-etal-2016,Zoppou-etal-2017} or a smooth approximation to them when the numerical method requires some smoothness of the solutions \cite{El-etal-2006,Mitsotakis-etal-2014,Mitsotakis-etal-2017}. Among these papers there are differences in the structures of the numerical solutions, with some demonstrating undulations in depth and velocity throughout the bore \cite{El-etal-2006,Zoppou-etal-2016,Zoppou-etal-2017} and others showing a constant depth and velocity state in the middle of the bore \cite{Hank-etal-2010-2034,Mitsotakis-etal-2014,Mitsotakis-etal-2017}.

The mean behaviour of numerical solutions to the dam-break problem for the Serre equations is consistent across the literature \cite{El-etal-2006,Zoppou-etal-2016,Zoppou-etal-2017,Hank-etal-2010-2034,Mitsotakis-etal-2014,Mitsotakis-etal-2017} and was demonstrated to be well approximated by the analytical solution to the dam-break problem by the shallow water wave equations \cite{Hank-etal-2010-2034,Mitsotakis-etal-2014}. Expressions for the leading wave amplitude and speed of an undular bore for the Serre equations were derived and verified for a range of undular bores by \citet{El-etal-2006}. These expressions were also shown to be valid for all the different structures found in the literature \cite{El-etal-2006,Mitsotakis-etal-2014}.

The first aim of this paper is to investigate and explain why different behaviour has been published in the literature for numerical solutions of the Serre equations for problems containing steep gradients. We find that the undulations of a bore can be damped to a constant depth and velocity state by the numerical diffusion introduced by the method, as is the case for \citet{Hank-etal-2010-2034}. Oscillation damping can also occur due to the particular smoothing of the initial conditions, as is the case for \citet{Mitsotakis-etal-2017}, \citet{El-etal-2006} and \citet{Mitsotakis-etal-2014}. We do find that over long time periods the Serre equations damp these oscillations as they propagate, but this natural decay is dominated by other factors in the literature.

The second aim of this paper is to assess the utility of the shallow water wave equations and the results of \citet{El-etal-2006} as guides for the evolution of an undular bore. We find that for a range of dam-break problems the analytical solution of the shallow water wave equations is a good approximation for the mean depth and velocity of the Serre equations, extending the findings of \citet{Hank-etal-2010-2034} and \citet{Mitsotakis-etal-2014} to a larger range of dam-break problems. It was also found that the results of \citet{El-etal-2006} are a good approximation to our numerical solutions.

The first aim of this paper is achieved by demonstrating that our numerical solutions are good approximations to the true solutions of the Serre equations. This is accomplished by demonstrating that as the resolution of a particular method is increased, the numerical solutions converge to a numerical solution with little error in the conservation of mass, momentum and energy. The numerical solution is also consistent across the five different numerical methods. Three of the methods are the first, second and third-order methods presented by \citet{Zoppou-etal-2017}. The first-order method is equivalent to the method of \citet{Hank-etal-2010-2034}. The fourth method is a recreation of the second-order method used by \citet{El-etal-2006}. Lastly, the fifth method is a second-order finite difference approximation to the Serre equations. 

The second aim is accomplished by comparing our verified numerical solutions to the analytical solutions of the shallow water wave equations and the Whitham modulation results presented by \citet{El-etal-2006}. 
 
The paper is organised as follows, in Section \ref{section:Serre Equations} the Serre equations and the quantities they conserve are presented. In Section \ref{section:smootheddambreak} the smoothed dam-break problem is defined, the measures of the relative difference between numerical solutions and the relative error in the conserved quantities are presented. The analytical solution of the shallow water wave equations and the expressions for the amplitude and speed of the leading wave of an undular bore are presented. In Section \ref{sec:nummeth} the numerical methods and their important properties are presented. In Section \ref{section:NumRes} the four different structures in the solutions of smoothed dam-break problem for the Serre equations are determined using verified numerical solutions. The verified numerical solutions are also used to evaluate how well the analytical solution of the shallow water wave equations captures the mean behaviour of the solution of the Serre equations for the dam-break problem. The Whitham modulations results are also compared to the verified numerical solutions to test their veracity.

\section{Serre Equations}
\label{section:Serre Equations}
The Serre equations can be derived by integrating the full inviscid incompressible Euler equations over the water depth \cite{Su-Gardener-1969-536}. They can also be derived as an asymptotic expansion of the Euler equations \cite{Bonneton-Lannes-2009-16601}. Assuming a constant horizontal bed, the one-dimensional Serre equations are \cite{Guyenne-etal-2014-169}
\begin{linenomath*}
\begin{subequations}\label{eq:Serre_nonconservative_form}
\begin{gather}
\dfrac{\partial h}{\partial t} + \dfrac{\partial (uh)}{\partial x} = 0
\label{eq:Serre_continuity}
\end{gather}
and
\begin{gather}
\underbrace{\underbrace{\dfrac{\partial (uh)}{\partial t} + \dfrac{\partial}{\partial x} \left ( u^2h + \dfrac{gh^2}{2}\right )}_{\text{Shallow Water Wave Equations}} + \underbrace{\dfrac{\partial}{\partial x} \left (  \dfrac{h^3}{3} \left [ \dfrac{\partial u }{\partial x} \dfrac{\partial u}{\partial x} - u\dfrac{\partial^2 u}{\partial x^2}  - \dfrac{\partial^2 u}{\partial x \partial t}\right ] \right )}_{\text{Dispersion Terms}} = 0.}_{\text{Serre Equations}}
\label{eq:Serre_momentum}
\end{gather}
\end{subequations}
\end{linenomath*}
Where $u(x,t)$ is the horizontal velocity over the depth of water $h(x,t)$, $g$ is the acceleration due to gravity, $x$ is the horizontal spatial variable and $t$ is time. 

The Serre equations are conservation laws for `mass' \eqref{eq:Serre_continuity}, `momentum' \eqref{eq:Serre_momentum} and the Hamiltonian \cite{Li-Y-2002,Green-Naghdi-1976-237}
\begin{linenomath*}
	\begin{gather}
	\label{eqn:Hamildef}
	\mathcal{H}(x,t) = \frac{1}{2} \left(hu^2 + \frac{h^3}{3} \left(\frac{\partial u}{\partial x}\right)^2 + gh^2\right)
	\end{gather}
\end{linenomath*}
which is the total energy.

The total amount of a quantity $q$ in a system in the spatial interval $[a,b]$ at a particular time $t$, is measured by
\begin{linenomath*}
\begin{gather*}
\label{eqn:Condef}
\mathcal{C}_q(t) = \int_{a}^{b} q(x,t)\, dx .
\end{gather*}
\end{linenomath*}
Conservation of a quantity $q$ implies that $\mathcal{C}_{q}(0) = \mathcal{C}_{q}(t)$ for all $t$ provided the interval is fixed and the system is closed. Our numerical methods should demonstrate conservation for the quantities $h$, $uh$ and $\mathcal{H}$.

\section{Smoothed Dam Break Problem}
\label{section:smootheddambreak}
In this section we define a class of initial condition problems, called the smoothed dam-break problem that we use throughout our numerical investigation. This class of initial conditions are used in the literature \cite{Mitsotakis-etal-2014,Mitsotakis-etal-2017} to smoothly approximate the discontinuous initial conditions of the dam-break problem, as some numerical methods require smoothness of the solutions.

\begin{figure}
\centering
\includegraphics[width=0.6\textwidth]{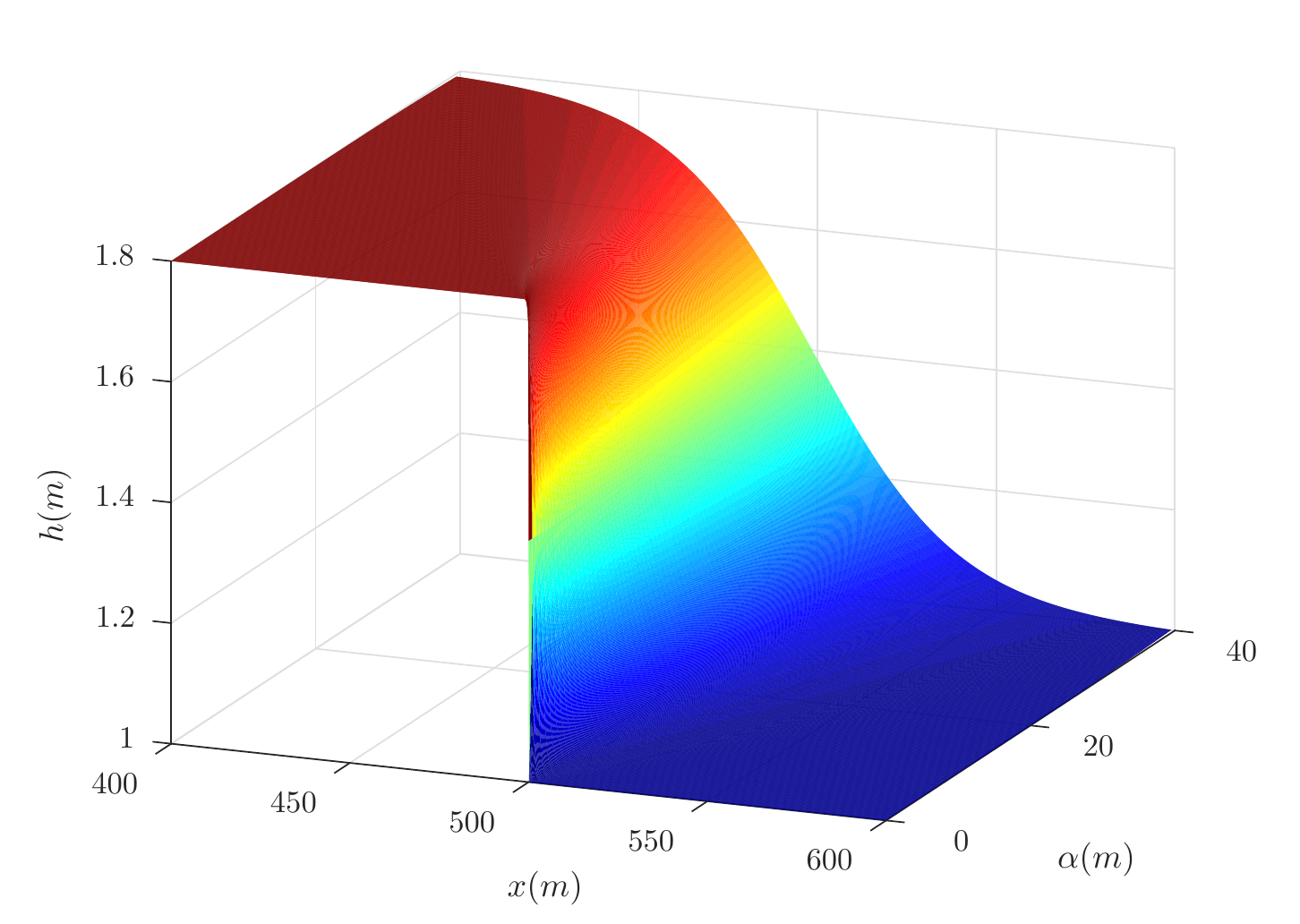}
\caption{Initial conditions for smooth dam-break problems with $h_0 = 1m$, $h_1 = 1.8m$, $x_0 =500m$ and various $\alpha$.}
\label{fig:dbsmoothinit}
\end{figure}

The smoothed dam-break problem has the following initial conditions
\begin{linenomath*}
\begin{subequations}
\begin{gather}
h(x,0) = h_0 + \frac{h_1 - h_0}{2}\left(1 + \tanh\left(\frac{x_0 - x}{\alpha}\right)\right)\; m,
\end{gather} 
and
\begin{gather}
u(x,0) = 0.0 \;m/s.
\end{gather}
\label{eq:sdbi}
\end{subequations}
\end{linenomath*}
This represents a smooth transition centred around $x_0$ between a water depth of $h_0$ on the right which is smaller than the water depth of $h_1$ on the left. Here $\alpha$ measures the distance over which approximately $46\%$ of that smooth transition between the two heights occurs. Decreasing $\alpha$ increases the steepness of the initial conditions as can be seen in Figure \ref{fig:dbsmoothinit} where $h_0=1m$ and $h_1=1.8m$. These are the same $h_0$ and $h_1$ values as those of the smoothed dam-break problem of \citet{El-etal-2006} and the dam-break problem of \citet{Hank-etal-2010-2034}.

There are no known analytical solutions of the Serre equations for the dam-break problem or an arbitrary smoothed dam-break problem. Therefore, to demonstrate that our numerical solutions converge we use the relative difference between numerical solutions. To demonstrate that our numerical solutions also have small errors in the conserved quantities we use the relative error of their conservation. Both of these measures are defined in this section.
\subsection{Assessing validity of Numerical Solutions}
To demonstrate that our numerical solutions converge to a solution with little error in the conserved quantities as the spatial resolution is increased we use two measures; the relative difference between numerical solutions of different resolutions and the error in the conservation of a quantity. The relative difference between numerical solutions measures their convergence, while the error in conservation measures how well the numerical solutions conserve the quantities $h$, $uh$ and $\mathcal{H}$.

We introduce the following notation for the spatial grids defined by $x_i$ and the temporal grids defined by $t^n$ upon which the numerical solutions are calculated. These grids are uniform so that $\Delta x = x_{i} - x_{i-1}$ for all $ i$ and $\Delta t = t^{n} - t^{n-1}$ for all $n$. We use subscripts and superscripts to denote where a quantity $q$ is evaluated in the following way $q_i^n = q(x_i,t^n)$. Finally, the $i$th cell is the interval [$x_i -\Delta x/2$,$x_i +\Delta x/2$] centred around $x_{i}$. 

\subsubsection{Convergence of Numerical Results}
To measure the convergence of the numerical solutions we ensured all grids had common locations to compare them by dividing $\Delta x$ by $2$ to create finer grids. Therefore, the finest grid with the smallest $\Delta x$ contains all the locations $x_i$ in any coarser grid. To measure the relative difference between quantities on these grids we compare them only on the coarser grid points $x_i$. For some quantity $q$ we have our numerical approximation to it on the finest grid $q^*$ and on the coarser grid $q'$, with the relative difference between the two being
\begin{linenomath*}
	\begin{gather}
	L_1^{q} = \dfrac{\sum_{i} \left| q'(x_i)  - q^*(x_i)\right|}{\sum_{i} \left| q^*(x_i)\right|}.
	\label{eq:L1def}
	\end{gather}
\end{linenomath*}

\subsubsection{Conserved Quantities}
To calculate the error in conservation of a quantity, we must first calculate the total amount of the conserved quantities for the initial conditions. For the smoothed dam-break problem the initial conditions \eqref{eq:sdbi} were integrated to obtain expressions for the total mass $C_{h}(0)$, the total momentum $C_{uh}(0)$ and the total Hamiltonian $C_{\mathcal{H}}(0)$. Provided $x_0$ is the midpoint of the spatial domain $\left[a,b \right]$ the total amounts for the conserved quantities are
\begin{linenomath*}
	\begin{subequations}
		\begin{gather*}
		\mathcal{C}_{h}(0) = \frac{h_1 + h_0}{2}\left(b- a\right),
		\label{eq:Chdef}
		\end{gather*}
		\begin{gather*}
		\mathcal{C}_{uh}(0) = 0
		\label{eq:Cuhdef}
		\end{gather*}
		and
		\begin{gather*}
		\mathcal{C}_{\mathcal{H}}(0) = \frac{g}{4} \left(h_0^2 - h_1^2 + \alpha\left(h_1 - h_0\right)^2\tanh\left(\frac{a - b}{2 \alpha}\right)\right).
		\label{eq:CHdef}
		\end{gather*}
		\label{eq:Canalyticvalues}	
	\end{subequations}
\end{linenomath*}

To calculate how well we approximate the total amount of a quantity $q$ in our numerical solution we fit a quartic interpolant of the primitive variables $h$ and $u$ over a cell utilising neighbouring cells and then apply Gaussian quadrature with 3 points. The amount of $q$ in each cell is summed across all cells to get the total amount of $q$ in our numerical solution at time $t$, which we call $\mathcal{C}^*_{q}(t)$. The error in conservation of a quantity $q$ for a numerical solution is
\begin{linenomath*}
	\begin{gather}
	C_1^q = \frac{\left| \mathcal{C}_{q}(0) - \mathcal{C}^*_{q}(t) \right| }{\left|\mathcal{C}_{q}(0)\right|}.
	\end{gather}
\end{linenomath*}
Note that for $uh$ the denominator is $0$ and that there is a flux of momentum due to the unequal heights at both ends of the domain. To resolve this issue for $uh$ the error in the conservation of $uh$ is measured by
\begin{linenomath*}
	\begin{gather}
	C_1^{uh} = \left| \mathcal{C}_{uh}(0) - \mathcal{C}^*_{uh}(t) - \frac{gt}{2}\left(h(b)^2 - h(a)^2\right)\right|  .
	\label{eq:C1def}
	\end{gather}
\end{linenomath*}

\subsection{Background for derived and observed comparisons}
It was demonstrated by \citet{Hank-etal-2010-2034} and \citet{Mitsotakis-etal-2014} that the analytical solution of the shallow water wave equations for the dam-break problem captures the mean behaviour of the numerical solutions of the Serre equations to the dam-break problem \cite{Hank-etal-2010-2034} and the smoothed dam-break problem \cite{Mitsotakis-etal-2014}. 

\citet{El-etal-2006} derived an expression for the long term amplitude of the leading wave of an undular bore $A^+$ for the Serre equations. Since the front of an undular bore decomposes into solitons, the speed of the leading wave $S^+$ can be calculated from its amplitude.

To be self contained we present the analytical solution of the shallow water wave equations to the dam-break problem and the expressions derived by \citet{El-etal-2006}.

\subsubsection{Shallow Water Wave Equation Analytical Solution}
For the dam-break problem the shallow water wave equations, which are the Serre equations with dispersive terms neglected, can be solved analytically. 

\begin{figure}
	\centering
	\includegraphics[width=0.5\textwidth]{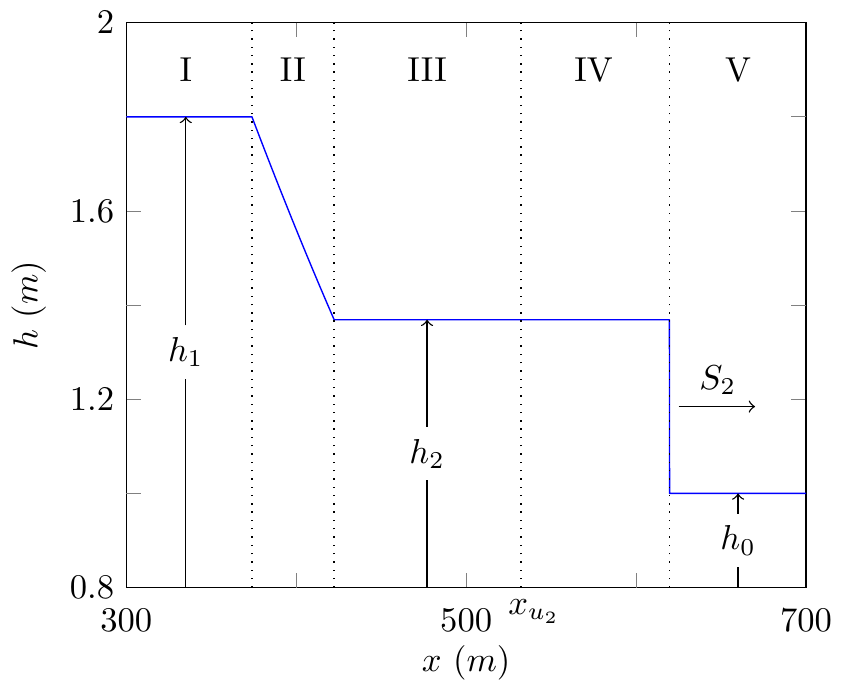}
	\caption{Analytical solution at $t=30s$ of the dam-break problem for the shallow water wave equations with $h_0 = 1m$, $h_1=1.8m$ and $x_0=500m$.}
	\label{fig:SWWEanadiagram}
\end{figure}

An example of the analytical solution of the shallow water wave equations for the dam-break problem is presented in Figure \ref{fig:SWWEanadiagram}. Region I is the undisturbed water upstream of the dam-break at constant height ($h_1$) and velocity ($0m/s$). Region II is the rarefaction fan connecting regions I and III. Regions III and IV represent the shock with constant height ($h_2$) and constant velocity ($u_2$), these regions are separated by $x_{u_2} = x_0 + u_2t$. Region V is the undisturbed water downstream at constant height ($h_0$) and velocity ($0m/s$) separated from Region IV by a shock which travels at velocity $S_2$. Expressions for the unknown quantities $h_2$, $u_2$ and $S_2$ in terms of $h_0$ and $h_1$ were given by \citet{Wu-etal-1999-1210} as
\begin{linenomath*}
\begin{subequations}
\begin{gather}
h_2 = \frac{h_0}{2} \left(\sqrt{1 + 8 \left(\frac{2h_2}{h_2 - h_0}\frac{\sqrt{h_1} - \sqrt{h_2}}{\sqrt{h_0}}\right)^2} - 1\right),
\end{gather}
	\begin{gather}
	u_2 = 2\left(\sqrt{gh_1} - \sqrt{gh_2}\right)
	\end{gather}
and
	\begin{gather}
	S_2 = \frac{h_2 u_2}{h_2 - h_0}.
	\end{gather}
\label{eq:WuSWWE}	
\end{subequations}
\end{linenomath*}
Applying \eqref{eq:WuSWWE} to our dam-break heights of interest; $h_0 =1m$ and $h_1 = 1.8m$ results in $h_2 = 1.36898m$ , $u_2 = 1.074975$ $m/s$ and $S_2 = 3.98835$ $m/s$ which are shown in Figure \ref{fig:SWWEanadiagram} for $t=30s$. The location of the front of the bore for the shallow water wave equations at time $t$ is thus $x_{S_2}(t) = x_0 + S_2 t$ so that $x_{S_2}(30s) = 619.6505m$.

\subsubsection{Whitham Modulation for Undular Bores of the Serre Equations}
\begin{figure}
	\centering
	\includegraphics[width=0.5\textwidth]{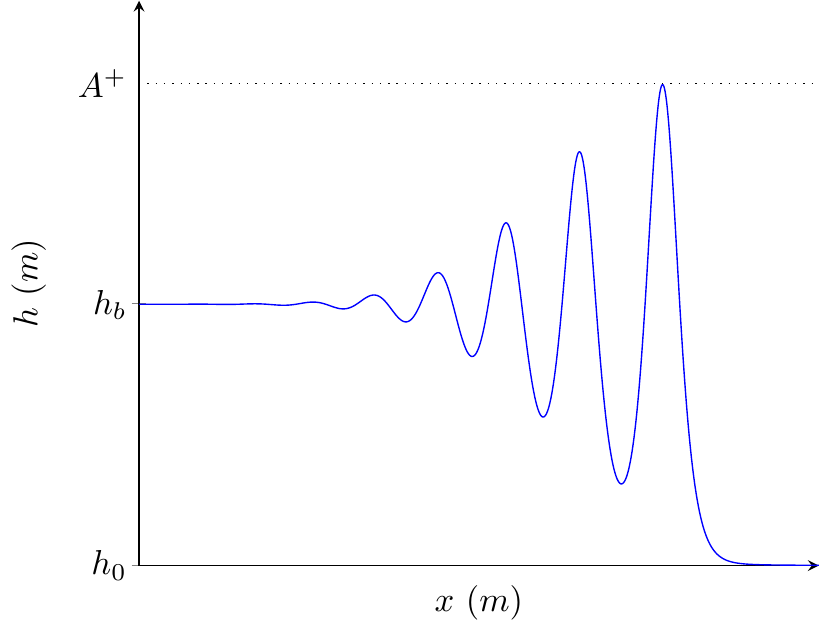}
	\caption{Demonstration of quantities obtained by Whitham modulation for undular bores of the Serre equations.}
	\label{fig:Serreanadiagram}
\end{figure}

Utilizing Whitham modulation theory for a one-phase periodic travelling wave an asymptotic analytical expression for the amplitude $A^+$ and speed $S^+$ of the leading wave was derived by \citet{El-etal-2006}. An example of an undular bore is shown in Figure \ref{fig:Serreanadiagram}.
The derived expressions for $A^+$ and $S^+$ are
\begin{linenomath*}
	\begin{subequations}
		\begin{gather}
		\frac{\Delta}{\left(A^+ + 1\right)^{1/4}} - \left(\frac{3}{4 -  \sqrt{A^+ + 1}}\right)^{21/10} \left(\frac{2}{1 + \sqrt{A^+ + 1}}\right)^{2/5} = 0
		\label{eq:aplusdef}
		\end{gather}
		and
		\begin{gather}
		S^+ = \sqrt{g \left(A^+ + 1\right)}
		\label{eq:splusdef}
		\end{gather}
		\label{eq:ELWhitMod}	
	\end{subequations}
\end{linenomath*}
where $\Delta = h_b / h_0$, and $h_b$ is the height of the bore. The height of the bore created by the dam-break problem in \eqref{eq:aplusdef} used by \citet{El-etal-2006} was
\begin{gather*}
\label{eqn:hrdef}
h_b = \frac{1}{4}\left(\sqrt{\frac{h_1}{h_0}} + 1\right)^2.
\end{gather*} 
For our dam-break heights of interest $h_0 = 1m$ and $h_1 = 1.8m$ we obtain $h_b = 1.37082m$, $\Delta = 1.37082$, $A^+ = 1.73998m$ and $S^+ = 4.13148m/s$. The location of the leading wave of an undular bore at time $t$ is then $x_{S^+}(t) = x_0 + S^+ t$ so that $x_{S^+}(30s) = 623.9444m$.

\section{Numerical Methods}
\label{sec:nummeth}
Five numerical schemes were used to investigate the behaviour of the Serre equations in the presence of steep gradients, the first ($\mathcal{V}_1$), second ($\mathcal{V}_2$) and third-order ($\mathcal{V}_3$) finite difference finite volume methods of  \citet{Zoppou-etal-2017}, the second-order finite difference method of \citet{El-etal-2006} ($\mathcal{E}$) and a second-order finite difference method ($\mathcal{D}$) that can be found in the Appendix.

The $\mathcal{V}_i$ methods are stable under a Courant-Friederichs-Lewy (CFL) condition  presented by \citet{Harten-etal-1983-357}. The $\mathcal{V}_i$ methods have demonstrated the appropriate order of convergence for smooth problems \cite{Zoppou-etal-2017}. Furthermore, $\mathcal{V}_2$ and $\mathcal{V}_3$ have been validated against experimental data containing steep gradients \cite{Zoppou-etal-2017}. The two methods $\mathcal{D}$ and $\mathcal{E}$ were found to be stable under the same CFL condition.

Generally, we found that $\mathcal{V}_1$ is the worst performing method due to its numerical diffusion \cite{Zoppou-etal-2017}. Of the high-order methods $\mathcal{E}$ is the worst performing, introducing dispersive errors.

\section{Numerical Results}
\label{section:NumRes}
We investigate the behaviour of the Serre equations in the presence of steep gradients by numerically solving the smoothed dam-break problem while varying the steepness of the initial conditions. As $\Delta x \rightarrow 0$ our numerical solutions should represent a good approximation of the true solution of the Serre equations. If our numerical solutions to a smoothed dam-break problem converge to the same numerical solution with little error in conservation of mass, momentum and energy as $\Delta x \rightarrow 0$ for each method, then this numerical solution is considered an accurate approximate solution to that smoothed dam-break problem for the Serre equations.

This process validates our numerical solutions for the smoothed dam-break problem, and thus validates our numerical methods to approximate the solution of the Serre equations in the presence of steep gradients, if it exists. With a validated model we can compare the numerical solution to the analytical solution of the shallow water wave equations for the dam-break problem and the results of \citet{El-etal-2006}.

Throughout most of this section we are interested in the numerical solution at $t=30s$ to the smoothed dam-break problem with $h_0 = 1m$, $h_1 = 1.8m$ and $x_0 = 500m$ while allowing for different $\alpha$ values. All numerical methods used $\Delta t = 0.01 \Delta x$ which is smaller than required by the CFL condition, ensuring stability of our schemes. The method $\mathcal{V}_2$ requires an input parameter to its slope limiter and this was chosen to be $\theta = 1.2$ \cite{Zoppou-etal-2017}. The spatial domain was $[0m,1000m]$ with the following Dirichlet boundary conditions, $u = 0m/s$ at both boundaries, $h =1.8m$ on the left and $h =1m$ on the right.

\subsection{Observed Structures of the Numerical Solutions}
\label{subsec:observedstructs}
\begin{figure}
	\centering
	\begin{subfigure}{0.5\textwidth}
		\includegraphics[width=\textwidth]{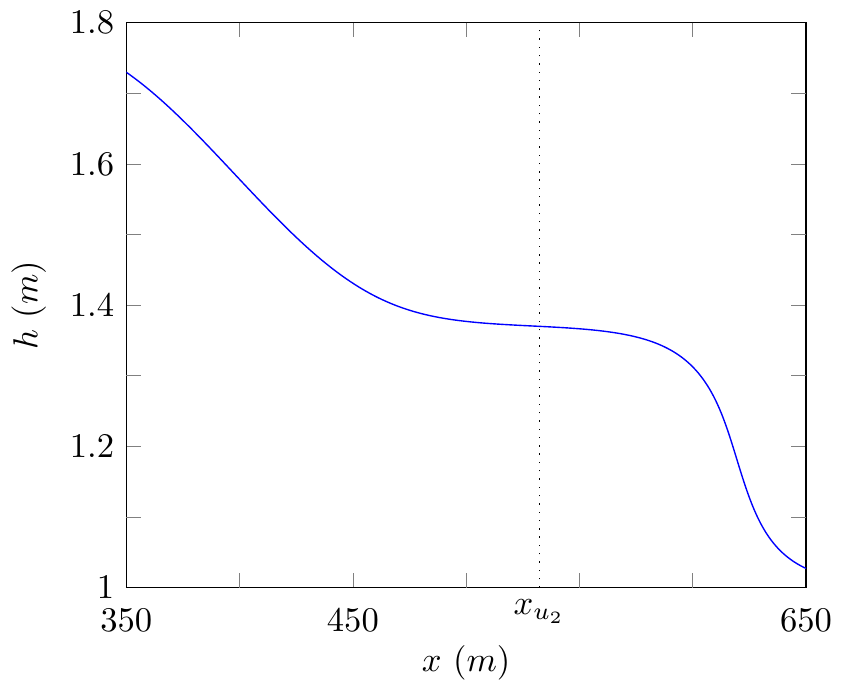}
		\subcaption*{\hspace{10 mm}$\mathcal{S}_1$  ($\alpha = 40m$)}
	\end{subfigure}%
	\begin{subfigure}{0.5\textwidth}
		\includegraphics[width=\textwidth]{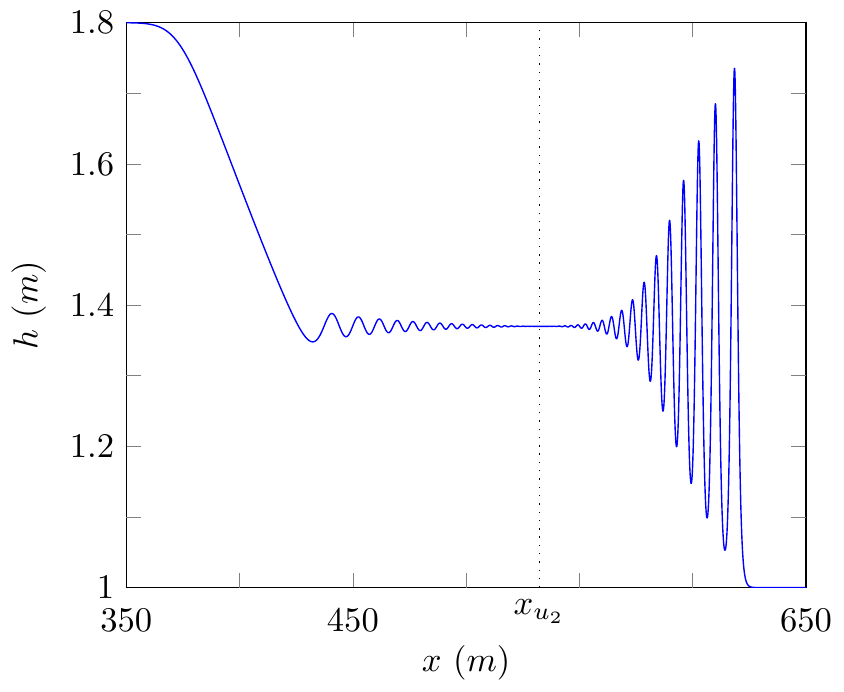}
		\subcaption*{\hspace{10 mm} $\mathcal{S}_2$ ($\alpha = 2m$)}
	\end{subfigure}
	\begin{subfigure}{0.5\textwidth}
		\includegraphics[width=\textwidth]{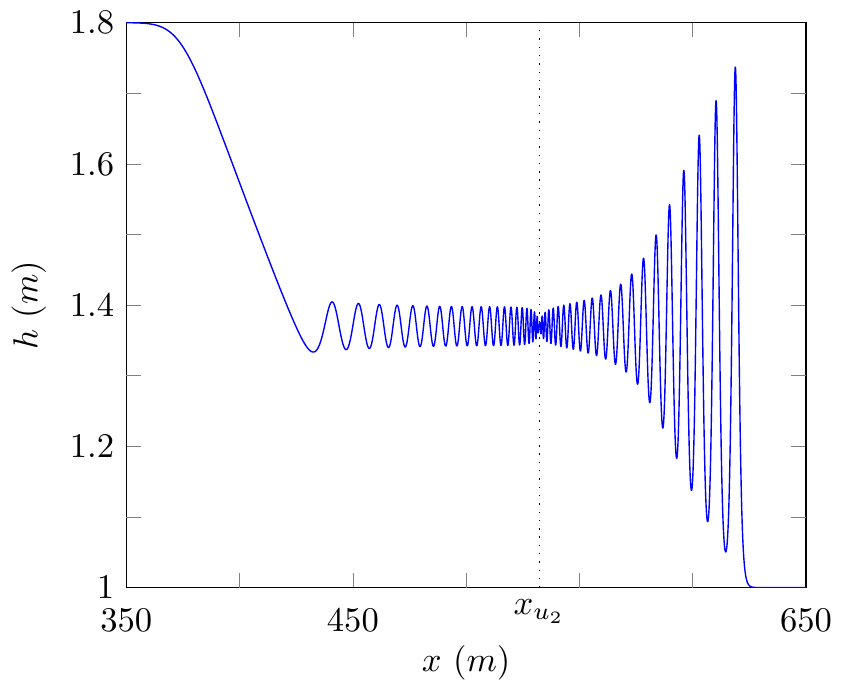}
		\subcaption*{\hspace{10 mm}$\mathcal{S}_3$ ($\alpha = 0.4m$)}
	\end{subfigure}%
	\begin{subfigure}{0.5\textwidth}
		\includegraphics[width=\textwidth]{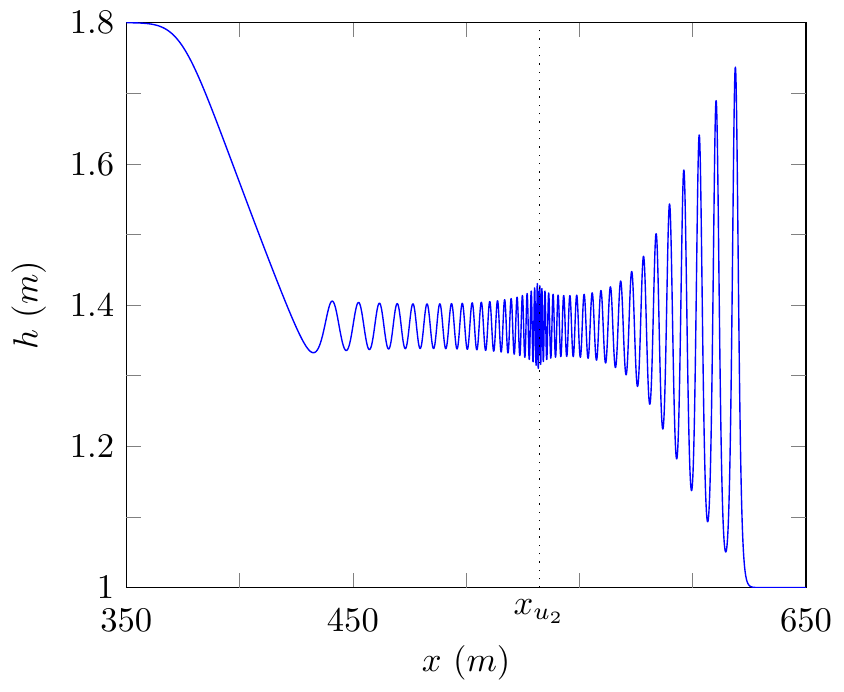}
		\subcaption*{\hspace{10 mm}$\mathcal{S}_4$ ($\alpha = 0.1m$)}
	\end{subfigure}	
	\caption{Numerical results of $\mathcal{V}_3$ with $\Delta x = 10/2^{11}m$ ({\color{blue} \solidrule}) at $t= 30s$ for various smooth dam-break problems demonstrating the different observed structures particularly around $x_{u_2}$ ({\color{black} \dotrule{4mm}}).}
	\label{fig:allstructs}
\end{figure}

We observe that there are four different structures for the converged to numerical solution depending on the chosen $\alpha$. They are the `\textit{non-oscillatory}' structure $\mathcal{S}_1$, the `\textit{flat}' structure $\mathcal{S}_2$, the `\textit{node}' structure $\mathcal{S}_3$ and the `\textit{growth}' structure $\mathcal{S}_4$. An example of each of these structures is shown in Figure \ref{fig:allstructs} which were obtained using $\mathcal{V}_3$ with $\Delta x = 10/2^{11}m$.

The four structures are identified by the dominant features of the numerical solutions in regions III and IV. They also correspond to different structures in the numerical solutions that have been presented in the literature. From Figure \ref{fig:allstructs} it can be seen that as $\alpha$ is decreased, steepening the initial conditions, the numerical solutions demonstrate an increase in the size and number of oscillations particularly around $x_{u_2}$. We observe that the difference between $\mathcal{S}_2$, $\mathcal{S}_3$ and $\mathcal{S}_4$ is the amplitude of the oscillations in regions III and IV.

\begin{figure}
	\centering
	\begin{subfigure}{0.49\textwidth}
		\includegraphics[width=\textwidth]{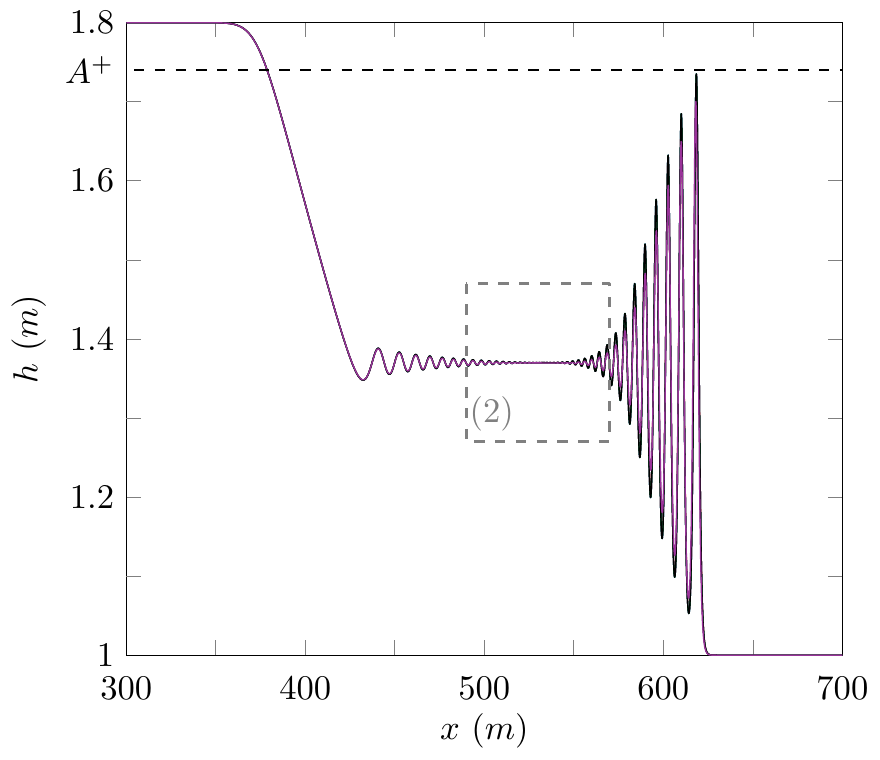}
	\end{subfigure}%
	\begin{subfigure}{0.5\textwidth}
		\includegraphics[width=\textwidth]{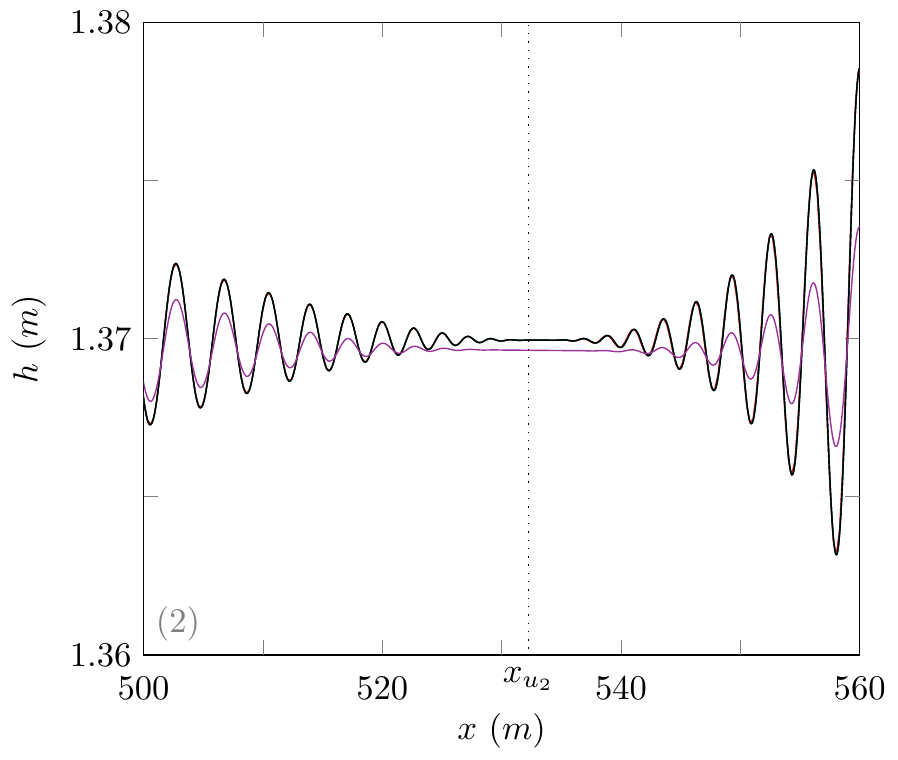}
	\end{subfigure}
	
	\caption{Numerical solutions of $\mathcal{D}$ ({\color{blue} \solidrule}), $\mathcal{E}$ ({\color{red} \solidrule}), $\mathcal{V}_3$ ({\color{green!60!black} \solidrule}), $\mathcal{V}_2$ ({\color{black} \solidrule}) and $\mathcal{V}_1$ ({\color{violet!60!white} \solidrule}) with $\Delta x = 10/2^{11}m$ at $t= 30s$ for the smooth dam-break problem with $\alpha =2m$. The Whitham modulation result for the leading wave height $A^+$ ({\color{black} \dashedrule}) and $x_{u_2}$ ({\color{black} \dotrule{4mm}}) are presented for comparison.}
	\label{fig:allmodels2struct}
\end{figure}

For the non-oscillatory and flat structures there is excellent agreement between all higher-order numerical methods at our highest resolution $\Delta x = 10/2^{11}m$. An illustration of this agreement is given in Figure \ref{fig:allmodels2struct} for $\mathcal{S}_2$ which is the most difficult to resolve of the two structures. However, the first-order method $\mathcal{V}_1$ suppresses oscillations present in the numerical solutions of other methods due to its diffusive errors \cite{Zoppou-etal-2017}. To resolve these oscillations with $\mathcal{V}_1$ much lower of values of $\Delta x$ are required.

\subsubsection{Non-oscillatory Structure}
The $\mathcal{S}_1$ ``non-oscillatory'' structure is the result of a large $\alpha$, which causes the front of this flow to not be steep enough to generate undulations over short time periods. As the system evolves the front will steepen due to non-linearity and undulations will develop.

\begin{figure}
	\centering
	\begin{subfigure}{0.5\textwidth}
		\includegraphics[width=\textwidth]{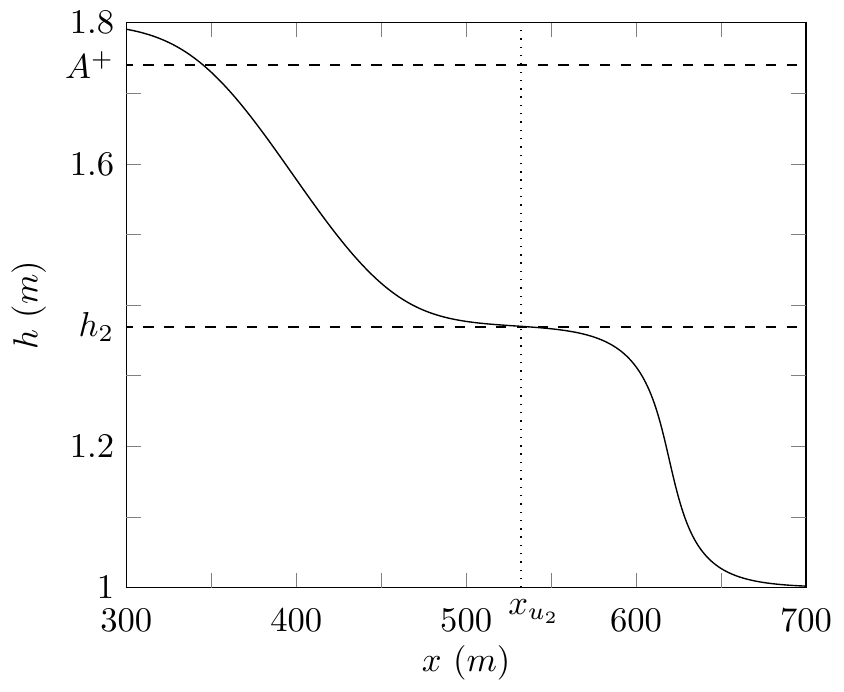}
	\end{subfigure}%
	
	\caption{Numerical solutions of $\mathcal{V}_3$ at $t= 30s$ for smooth dam-break problem with $\alpha = 40m$ for $\Delta x = 10/2^{10}m$ ({\color{blue} \solidrule}), $10/2^8m$ ({\color{red} \solidrule}), $10/2^6m$ ({\color{green!60!black} \solidrule}) and $10/2^{4}m$ ({\color{black} \solidrule}). The important quantities $A^+$ ({\color{black} \dashedrule}), $h_2$ ({\color{black} \dashedrule}) and $x_{u_2}$ ({\color{black} \dotrule{4mm}}) are also presented.}
	\label{fig:o3a1dxlimnonexp}
\end{figure}

The structure $\mathcal{S}_1$ is not present in the literature as no authors chose large enough $\alpha$ because, such a large $\alpha$ poorly approximates the dam-break problem. An example of this structure can be seen in Figure \ref{fig:o3a1dxlimnonexp} for $\alpha = 40m$ using $\mathcal{V}_3$ with various $\Delta x$ values. Because this is not a very steep problem all numerical results are visually identical for all $\Delta x < 10 / 2^4m$.

From Table \ref{tab:L1C1} it can be seen that not only have these solutions converged visually but the $L_1$ measures demonstrate that we have reached convergence to round-off error by $\Delta x = 10 / 2^8m$ after which the relative difference between numerical solutions plateau. 

Table \ref{tab:L1C1} also demonstrates that the error in conservation of the numerical solutions are at round-off error for $h$ and $\mathcal{H}$. The conservation of $uh$ is poor because the smoothed dam-break has such a large $\alpha$ that $h(0m) \neq 1.8m$ and $h(1000m) \neq 1m$, causing unequal fluxes in momentum at the boundaries. 

As stated above when $\Delta x = 10/2^{11}m$ the numerical solutions from all methods are identical for this smoothed dam-break problem. 

The convergence of the numerical solutions as $\Delta x \rightarrow 0$ to a numerical solution with small error in conservation, independent of the method demonstrates that we have accurately solved the smoothed dam-break problem with $\alpha = 40m$. Therefore, the $\mathcal{S}_1$ structure should be observed in the solutions of the Serre equations for the smoothed dam-break problem for sufficiently large $\alpha$. 

\begin{table}
\begin {tabular}{c c c c c c c}%
$\alpha$& $\Delta x$&$C_1^{h}$&$C_1^{uh}$&$C_1^{\mathcal {H}}$&$L_1^{h}$&$L_1^{u}$\\%
\hline\hline\\
$40$ &$10/2^{4}$& \pgfutilensuremath {2.00\cdot 10^{-11}}&\pgfutilensuremath {1.77\cdot 10^{-6}}&\pgfutilensuremath {1.23\cdot 10^{-8}}&\pgfutilensuremath {1.74\cdot 10^{-7}}&\pgfutilensuremath {2.90\cdot 10^{-6}}\\%
$40$ &$10/2^{6}$& \pgfutilensuremath {1.07\cdot 10^{-11}}&\pgfutilensuremath {1.50\cdot 10^{-6}}&\pgfutilensuremath {1.49\cdot 10^{-10}}&\pgfutilensuremath {2.57\cdot 10^{-9}}&\pgfutilensuremath {4.19\cdot 10^{-8}}\\%
$40$ &$10/2^{8}$& \pgfutilensuremath {8.77\cdot 10^{-13}}&\pgfutilensuremath {5.49\cdot 10^{-7}}&\pgfutilensuremath {3.77\cdot 10^{-13}}&\pgfutilensuremath {6.08\cdot 10^{-11}}&\pgfutilensuremath {5.28\cdot 10^{-10}}\\%
$40$ &$10/2^{10}$& \pgfutilensuremath {1.77\cdot 10^{-11}}&\pgfutilensuremath {2.21\cdot 10^{-8}}&\pgfutilensuremath {3.56\cdot 10^{-11}}&\pgfutilensuremath {2.54\cdot 10^{-11}}&\pgfutilensuremath {6.49\cdot 10^{-11}}\\%
\\ \hline \\
$2$ &$10/2^{4}$&\pgfutilensuremath {4.90\cdot 10^{-14}}&\pgfutilensuremath {5.10\cdot 10^{-3}}&\pgfutilensuremath {8.69\cdot 10^{-4}}&\pgfutilensuremath {5.02\cdot 10^{-3}}&\pgfutilensuremath {6.77\cdot 10^{-2}}\\%
$2$ &$10/2^{6}$&\pgfutilensuremath {2.51\cdot 10^{-13}}&\pgfutilensuremath {2.18\cdot 10^{-4}}&\pgfutilensuremath {6.58\cdot 10^{-5}}&\pgfutilensuremath {4.14\cdot 10^{-4}}&\pgfutilensuremath {5.20\cdot 10^{-3}}\\%
$2$ &$10/2^{8}$&\pgfutilensuremath {9.81\cdot 10^{-13}}&\pgfutilensuremath {7.72\cdot 10^{-7}}&\pgfutilensuremath {5.01\cdot 10^{-7}}&\pgfutilensuremath {6.00\cdot 10^{-6}}&\pgfutilensuremath {7.59\cdot 10^{-5}}\\%
$2$ &$10/2^{10}$&\pgfutilensuremath {3.95\cdot 10^{-12}}&\pgfutilensuremath {5.56\cdot 10^{-9}}&\pgfutilensuremath {6.13\cdot 10^{-9}}&\pgfutilensuremath {1.76\cdot 10^{-7}}&\pgfutilensuremath {2.33\cdot 10^{-6}}\\
 \\ \hline \\
$0.4$ &$10/2^{4}$&\pgfutilensuremath {9.00\cdot 10^{-14}}&\pgfutilensuremath {4.82\cdot 10^{-3}}&\pgfutilensuremath {1.02\cdot 10^{-3}}&\pgfutilensuremath {6.79\cdot 10^{-3}} $\dagger$&\pgfutilensuremath {9.93\cdot 10^{-2}} $\dagger$\\%
$0.4$ &$10/2^{6}$&\pgfutilensuremath {2.40\cdot 10^{-13}}&\pgfutilensuremath {2.41\cdot 10^{-4}}&\pgfutilensuremath {1.11\cdot 10^{-4}}&\pgfutilensuremath {8.89\cdot 10^{-4}} $\dagger$&\pgfutilensuremath {1.13\cdot 10^{-2}} $\dagger$\\%
$0.4$ &$10/2^{8}$&\pgfutilensuremath {9.68\cdot 10^{-13}}&\pgfutilensuremath {7.57\cdot 10^{-7}}&\pgfutilensuremath {2.25\cdot 10^{-6}}&\pgfutilensuremath {1.53\cdot 10^{-5}} $\dagger$&\pgfutilensuremath {1.91\cdot 10^{-4}} $\dagger$\\%
$0.4$ &$10/2^{10}$&\pgfutilensuremath {3.91\cdot 10^{-12}}&\pgfutilensuremath {4.95\cdot 10^{-9}}&\pgfutilensuremath {2.01\cdot 10^{-8}}&\pgfutilensuremath {3.61\cdot 10^{-7}} $\dagger$&\pgfutilensuremath {5.00\cdot 10^{-6}} $\dagger$\\%
\\ \hline \\
$0.1$ &$10/2^{4}$&\pgfutilensuremath {7.60\cdot 10^{-14}}&\pgfutilensuremath {4.82\cdot 10^{-3}}&\pgfutilensuremath {1.06\cdot 10^{-3}}&\pgfutilensuremath {7.04\cdot 10^{-3}} $\dagger$&\pgfutilensuremath {1.02\cdot 10^{-1}} $\dagger$\\%
$0.1$ &$10/2^{6}$&\pgfutilensuremath {2.40\cdot 10^{-13}}&\pgfutilensuremath {2.39\cdot 10^{-4}}&\pgfutilensuremath {1.44\cdot 10^{-4}}&\pgfutilensuremath {1.02\cdot 10^{-3}} $\dagger$&\pgfutilensuremath {1.28\cdot 10^{-2}} $\dagger$\\%
$0.1$ &$10/2^{8}$&\pgfutilensuremath {9.79\cdot 10^{-13}}&\pgfutilensuremath {2.21\cdot 10^{-7}}&\pgfutilensuremath {1.20\cdot 10^{-5}}&\pgfutilensuremath {2.86\cdot 10^{-5}} $\dagger$&\pgfutilensuremath {3.46\cdot 10^{-4}} $\dagger$\\%
$0.1$ &$10/2^{10}$&\pgfutilensuremath {3.92\cdot 10^{-12}}&\pgfutilensuremath {4.46\cdot 10^{-8}}&\pgfutilensuremath {7.61\cdot 10^{-7}}&\pgfutilensuremath {4.99\cdot 10^{-7}} $\dagger$&\pgfutilensuremath {6.40\cdot 10^{-6}} $\dagger$\\ \\ \hline\\%
\end {tabular}%
\caption{All errors in conservation $C^{q}_1$ for the conserved quantities and relative differences $L^{q}_1$ of the primitive variables for numerical solutions of $\mathcal{V}_3$. $L^{q}_1$ uses the numerical solution with $\Delta x = 10/2^{11}m$ as the high resolution basis of comparison and $\dagger$ indicates the omission of the interval [$520m$, $540m$] from the comparison.}
\label{tab:L1C1}
\end{table}

\subsubsection{Flat Structure}
\begin{figure}
		\begin{subfigure}{0.49\textwidth}
			\includegraphics[width=\textwidth]{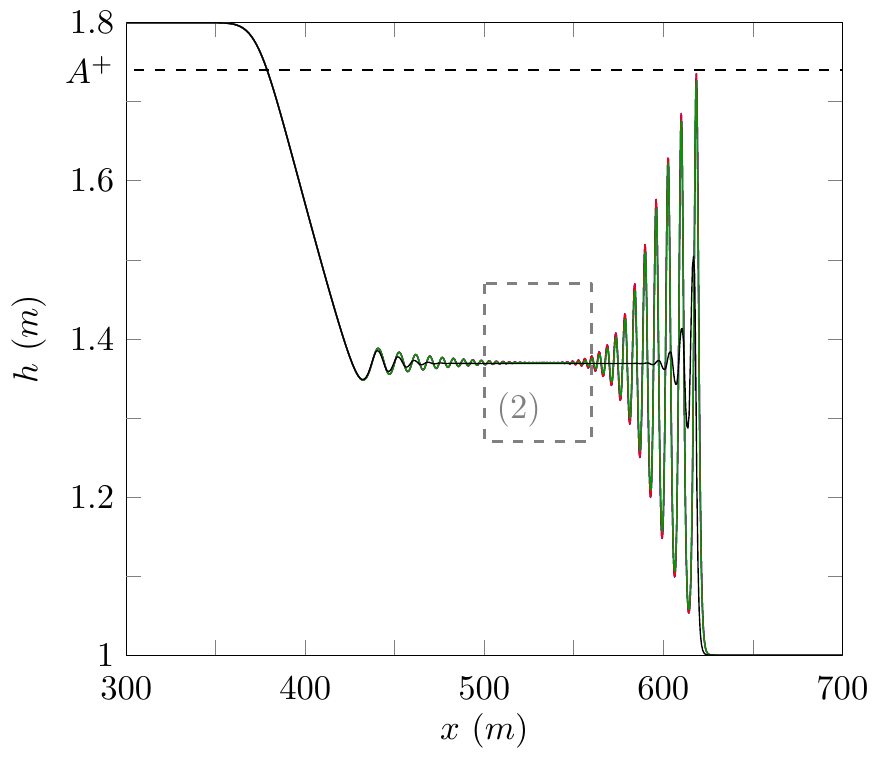}
		\end{subfigure}%
		\begin{subfigure}{0.5\textwidth}
			\includegraphics[width=\textwidth]{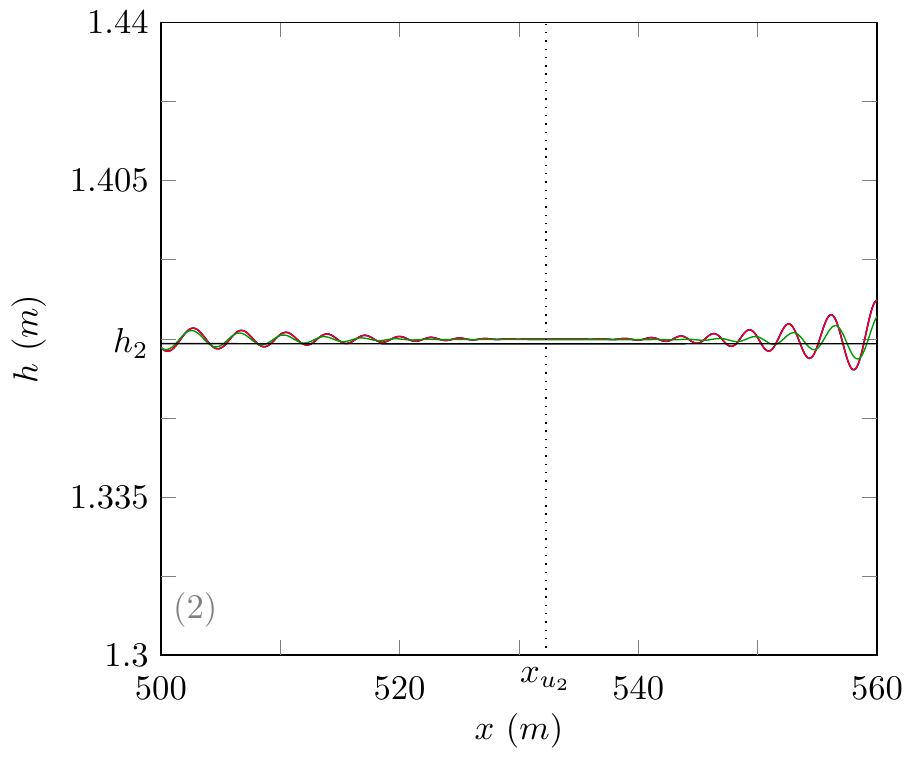}
		\end{subfigure}
\caption{Numerical solutions of $\mathcal{V}_3$ at $t= 30s$ for the smooth dam-break problem with $\alpha = 2m$ for $\Delta x = 10/2^{10}m$ ({\color{blue} \solidrule}), $10/2^8m$ ({\color{red} \solidrule}), $10/2^6m$ ({\color{green!60!black} \solidrule}) and $10/2^{4}m$ ({\color{black} \solidrule}). The important quantities $A^+$ ({\color{black} \dashedrule}), $h_2$ ({\color{black} \dashedrule}) and $x_{u_2}$ ({\color{black} \dotrule{4mm}}) are also presented.}
\label{fig:o3a6dxlimflatexp}
\end{figure}

The most common structure observed in the literature \cite{Hank-etal-2010-2034,Mitsotakis-etal-2014,Mitsotakis-etal-2017} is the ``flat structure'' $\mathcal{S}_2$. It is observed when the initial conditions are steep enough such that the bore that develops has undulations. This structure consists of oscillations in regions III and IV which are separated by a constant height state around $x_{u_2}$. An example of the $\mathcal{S}_2$ structure can be seen in the numerical solutions presented in Figure \ref{fig:o3a6dxlimflatexp} where $\alpha = 2m$.

As $\Delta x$ decreases the numerical solutions converge so that by $\Delta x = 10 / 2^8m$ the solutions for higher $\Delta x$ are visually identical. Table \ref{tab:L1C1} demonstrates that although we have convergence visually, the $L_1$ measures are still decreasing and are larger than round-off error. Likewise the $C_1$ measures are still decreasing and have only reached round-off error for $h$. This indicates that to attain full convergence of the numerical solutions of this smoothed dam-break problem down to round-off error using $\mathcal{V}_3$ would require an even smaller $\Delta x$. The relative difference between numerical solutions is small and the numerical solutions exhibit good conservation. Therefore, our highest resolution numerical solution is a good approximation to any numerical solutions with lower $\Delta x$ values. Figure \ref{fig:allmodels2struct} demonstrates that at $\Delta x = 10 / 2^{11}m$ the numerical solutions of all higher order methods are the same.

These results demonstrate that our highest resolution numerical solution is an accurate approximate solution of the Serre equations for the smoothed dam-break problem with $\alpha = 2m$. This implies that the $\mathcal{S}_2$ structure should be observed in solutions of the Serre equations for smooth dam-break problems with similar $\alpha$ values.

These numerical solutions compare well with those of \citet{Mitsotakis-etal-2014} who use the same $\alpha$ but different $h_0$ and $h_1$ values and observe the $\mathcal{S}_2$ structure. We found that we observed this structure for all numerical method's numerical solutions to the smoothed dam-break problem with $\alpha$ values as low as $1m$ and $\Delta x = 10/2^{11}m$. The numerical solutions of \citet{Mitsotakis-etal-2017} use $\alpha=1m$ but different heights and observe the structure $\mathcal{S}_2$. Therefore \citet{Mitsotakis-etal-2017} and \citet{Mitsotakis-etal-2014} observe the $\mathcal{S}_2$ structure in their numerical results due to their choice of $\alpha$ for the smoothed dam-break problem. 

\begin{figure}
	\centering
	\begin{subfigure}{\textwidth}
		\centering
		\includegraphics[width=0.5\textwidth]{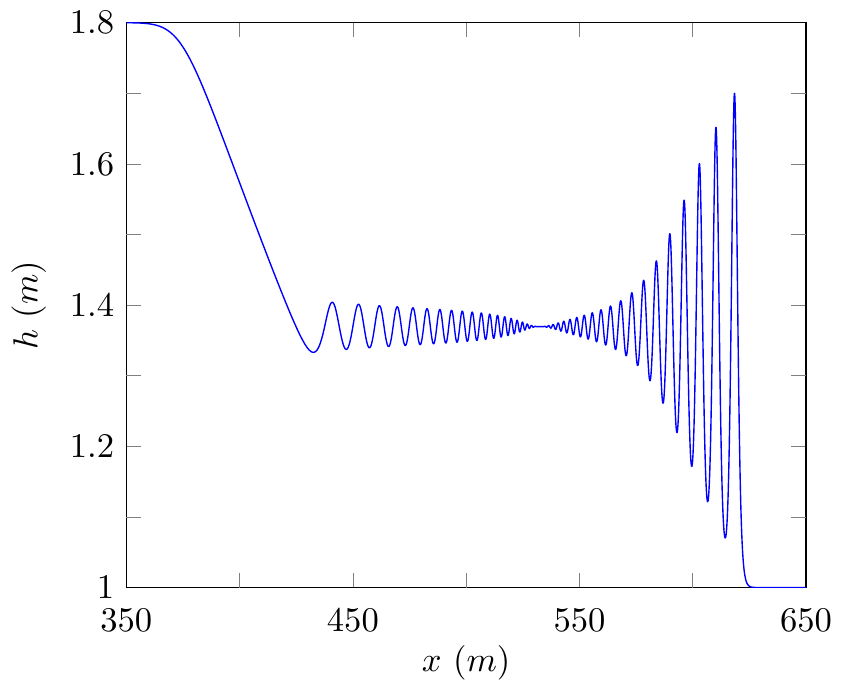}
	\end{subfigure}
	\caption{Numerical solution of $\mathcal{V}_1$ at $t= 30s$ for the smooth dam-break problem with $\alpha = 0.001m$ for $\Delta x = 10/2^{11}m$ ({\color{blue} \solidrule}).}
	\label{fig:o1highdchigha}
\end{figure}

The first-order method $\mathcal{V}_1$ is diffusive \cite{Zoppou-etal-2017} and damps oscillations that are present in the numerical solutions of higher-order methods as in Figure \ref{fig:allmodels2struct}. We find that for any smoothed dam-break problem with $\alpha \le 4m$ and the dam-break problem only the $\mathcal{S}_2$ structure is observed for the numerical solutions of $\mathcal{V}_1$ at $t=30s$ with $\Delta x = 10/2 ^{11}m$. This is evident in Figure \ref{fig:o1highdchigha} with the numerical solutions of $\mathcal{V}_1$ using our finest grid where $\Delta x = 10/2^{11}m$ on our steepest initial conditions where $\alpha = 0.001m$. Therefore, \citet{Hank-etal-2010-2034} using the diffusive $\mathcal{V}_1$ with their chosen $\Delta x$ and $\Delta t$, which are larger than our $\Delta x$ and $\Delta t$ could only observe the $\mathcal{S}_2$ structure. 

\subsubsection{Node Structure}
\begin{figure}
\centering
\begin{subfigure}{\textwidth}
	\centering
	\includegraphics[width=0.5\textwidth]{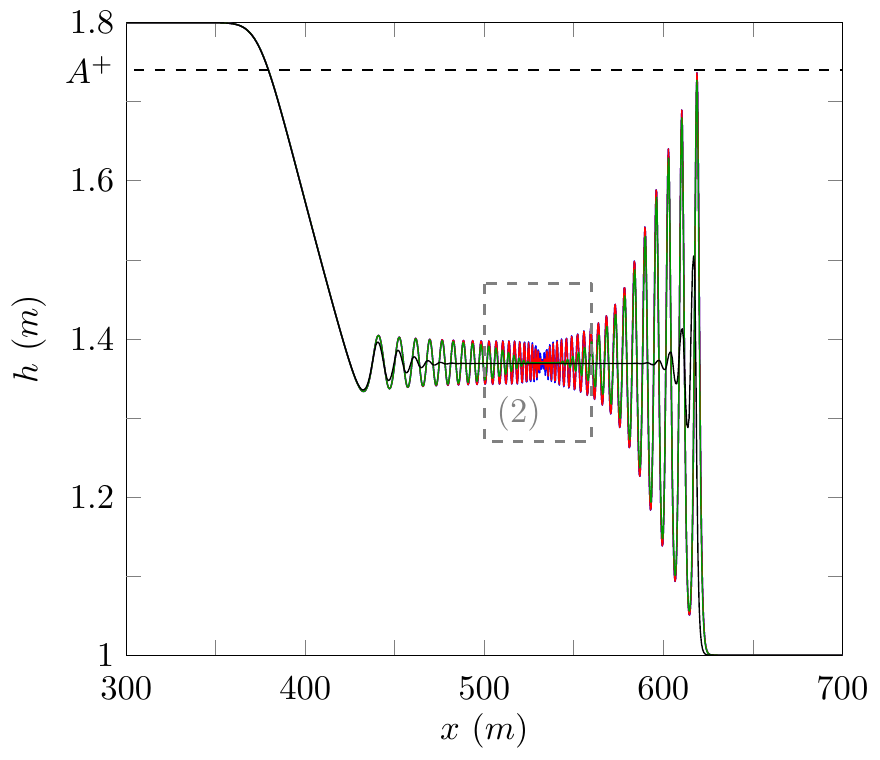}
\end{subfigure}
\begin{subfigure}{\textwidth}
	\includegraphics[width=0.5\textwidth]{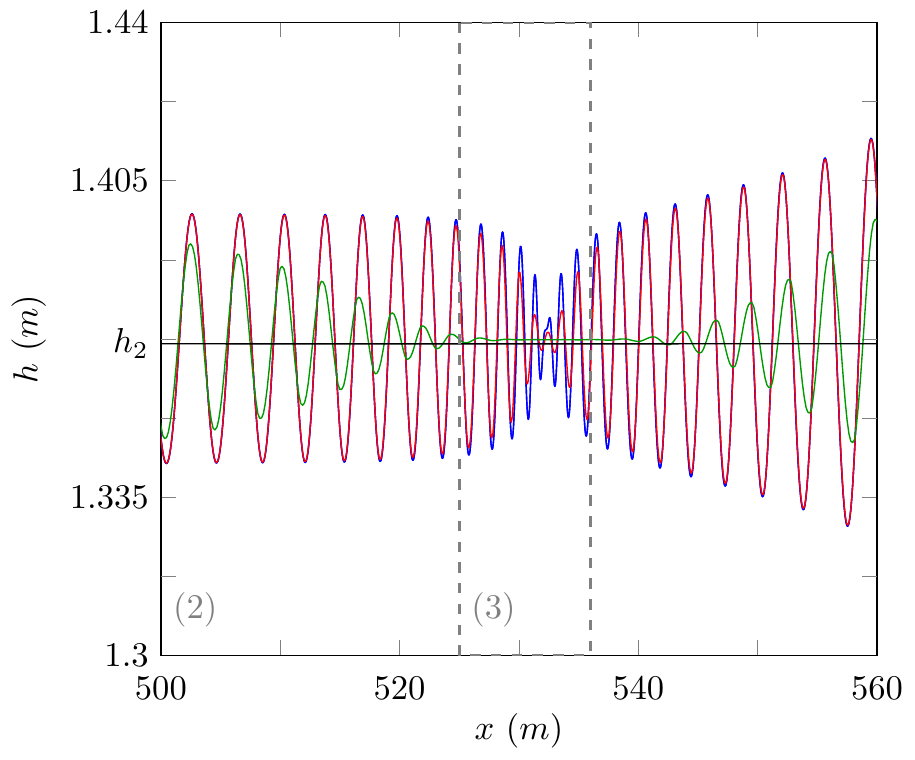}
	\includegraphics[width=0.5\textwidth]{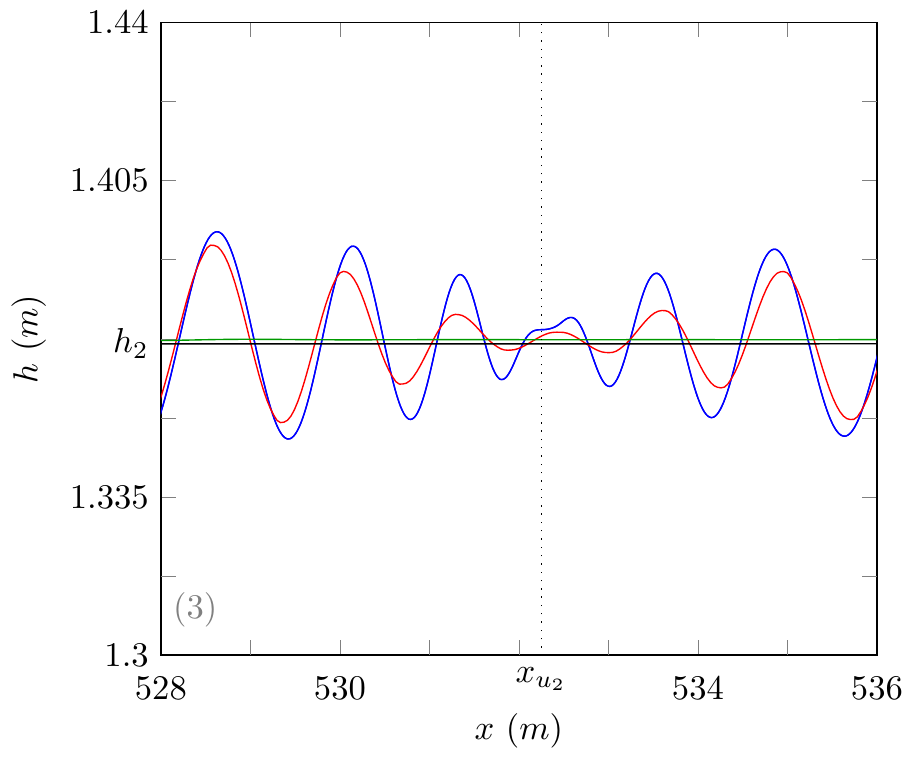}
\end{subfigure}
\caption{Numerical solutions of $\mathcal{V}_3$ at $t= 30s$ for the smooth dam-break problem with $\alpha = 0.4m$ for $\Delta x = 10/2^{10}m$ ({\color{blue} \solidrule}), $10/2^8m$ ({\color{red} \solidrule}), $10/2^6m$ ({\color{green!60!black} \solidrule}) and $10/2^{4}m$ ({\color{black} \solidrule}). The important quantities $A^+$ ({\color{black} \dashedrule}), $h_2$ ({\color{black} \dashedrule}) and $x_{u_2}$ ({\color{black} \dotrule{4mm}}) are also presented.}
\label{fig:o3a9dxlimcdexp}
\end{figure}

The ``node'' structure, $\mathcal{S}_3$ was observed by \citet{El-etal-2006}. The $\mathcal{S}_3$ structure has oscillations throughout regions III and IV that decay to a node at $x_{u_2}$ as can be seen in Figure \ref{fig:o3a9dxlimcdexp} where $\alpha = 0.4m$.

Figure \ref{fig:o3a9dxlimcdexp} demonstrates that our numerical solutions have not converged, however this is only in the area around $x_{u_2}$. Due to the large difference in numerical solutions around $x_{u_2}$ the $L_1$ measure over the area around $x_{u_2}$ would not be insightful. However, by omitting this region we can gain some knowledge about how well our solutions agree away from $x_{u_2}$. This was performed for the relevant $L_1$ measures in Table \ref{tab:L1C1} by omitting the interval [$520m$, $540m$]. These modified $L_1$ measures demonstrate that while our numerical results have visually converged outside this interval, they have not converged down to round-off error. 

Table \ref{tab:L1C1} demonstrates that the $C_1$ measures are still decreasing and have only attained round-off error for $h$. Therefore, to resolve the desired convergence of the numerical solutions to one with small error in conservation using $\mathcal{V}_3$ would require even smaller $\Delta x$ values.

\begin{figure}
	\centering
	\begin{subfigure}{\textwidth}
		\hspace{1mm}
		\includegraphics[width=0.49\textwidth]{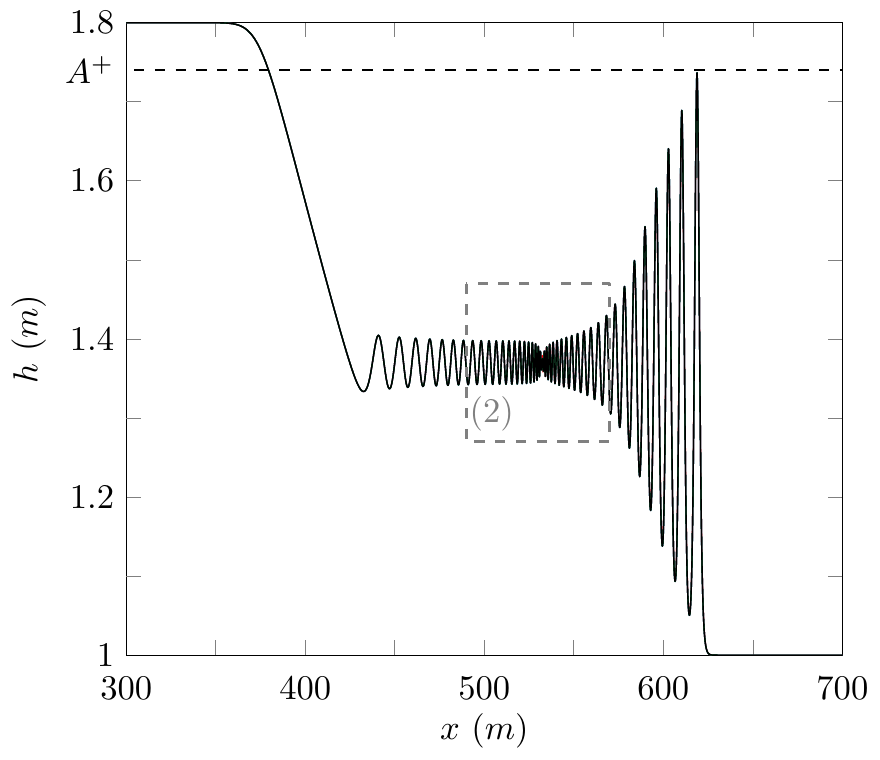}
		\includegraphics[width=0.5\textwidth]{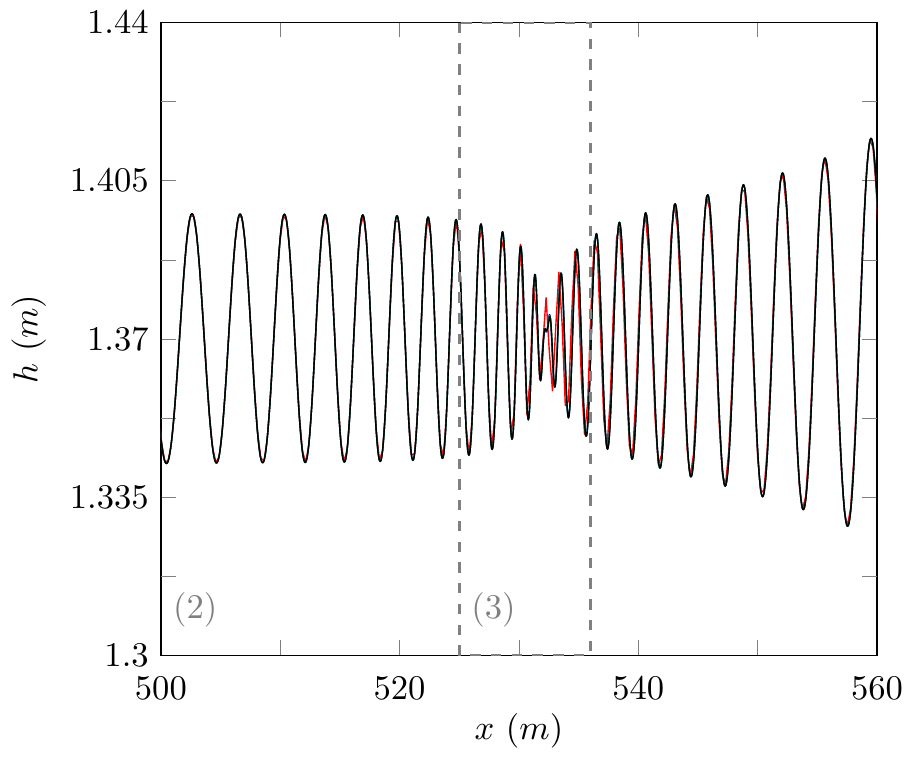}
	\end{subfigure}
	\begin{subfigure}{\textwidth}
		\includegraphics[width=0.5\textwidth]{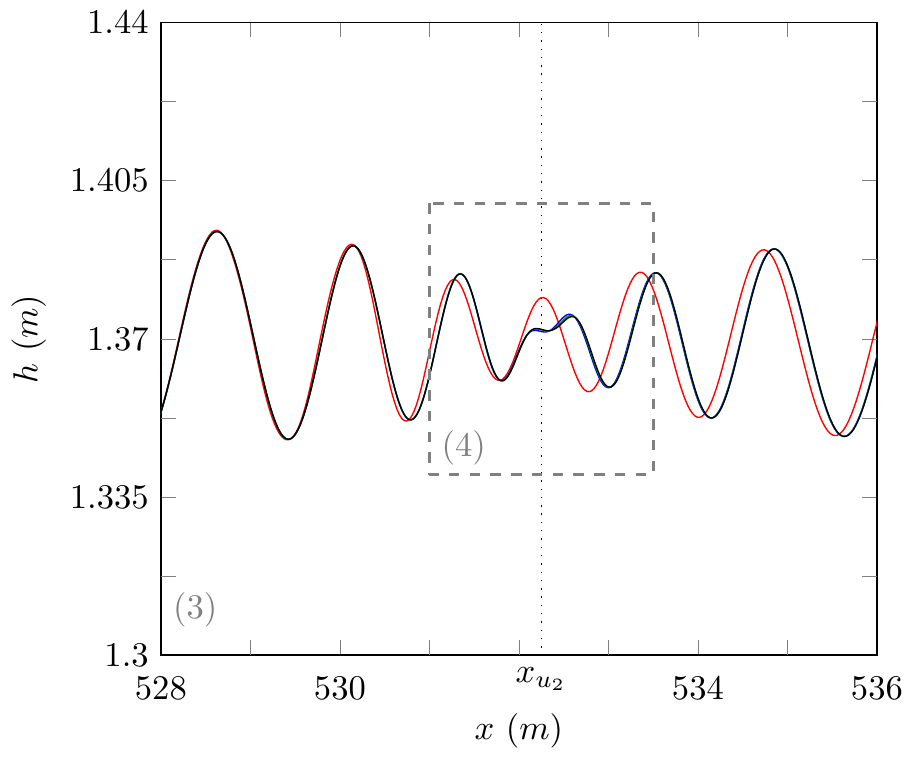}
		\includegraphics[width=0.49\textwidth]{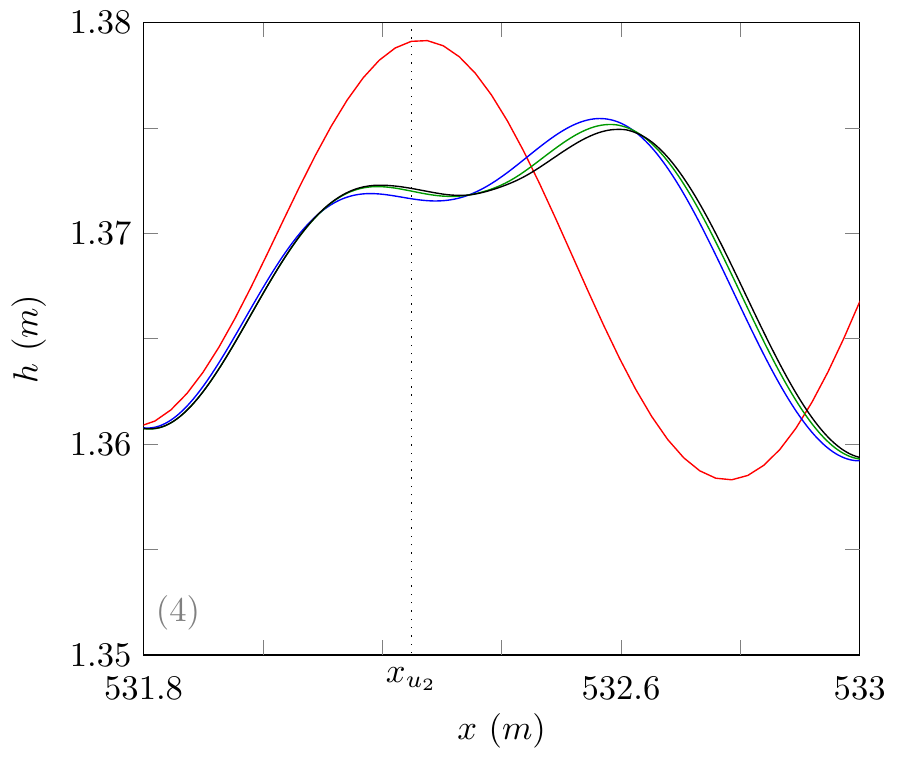}
	\end{subfigure}
	\caption{Numerical solutions of $\mathcal{D}$ ({\color{blue} \solidrule}), $\mathcal{E}$ ({\color{red} \solidrule}), $\mathcal{V}_3$ ({\color{green!60!black} \solidrule}) and $\mathcal{V}_2$ ({\color{black} \solidrule}) at $t=30s$ with $\Delta x = 10/2^{11}m$ for the smoothed dam-break problem with $\alpha = 0.4m$. The important quantities $A^+$ ({\color{black} \dashedrule}) and $x_{u_2}$ ({\color{black} \dotrule{4mm}}) are also presented.}
	\label{fig:o3a9allmodels}
\end{figure}

There is good agreement across different numerical methods for $\Delta x = 10/2^{11}m$ as can be seen in Figure \ref{fig:o3a9allmodels}. In particular all the higher-order methods exhibit the same structure and only disagree in a very small region around $x_{u_2}$. We observe that the numerical solution of the worst higher-order method $\mathcal{E}$ has not converged well to the numerical solutions of the other higher-order methods. 

We have only obtained a good approximation to the desired numerical solution as $\Delta x \rightarrow 0$ away from $x_{u_2}$. However, our highest resolution numerical solutions from various higher-order methods are very similar. This suggests that again although we do not have full convergence, our highest resolution numerical solution is a good approximation to the desired numerical solution over the whole domain. Therefore, our highest resolution numerical solutions are an accurate representation of the solutions of the Serre equations for this smoothed dam-break problem. Therefore, the $\mathcal{S}_3$ structure should be observed in the solutions of the Serre equations for the smoothed dam-break problem with $\alpha = 0.4m$.

These numerical solutions support the findings of \citet{El-etal-2006} who also use some smoothing \cite{El-Hoefer-2016-11} but do not report what smoothing was performed. Using their method $\mathcal{E}$ and similar $\Delta x$ to \citet{El-etal-2006} we observe the $\mathcal{S}_4$ ``growth'' structure in the numerical solution for $\alpha$ values smaller than $0.1m$, indicating that the smoothing performed by \citet{El-etal-2006} limited their observed behaviour to just the $\mathcal{S}_3$ structure. 

\subsubsection{Growth Structure}
\begin{figure}
\centering
\begin{subfigure}{\textwidth}
	\centering
	\includegraphics[width=0.5\textwidth]{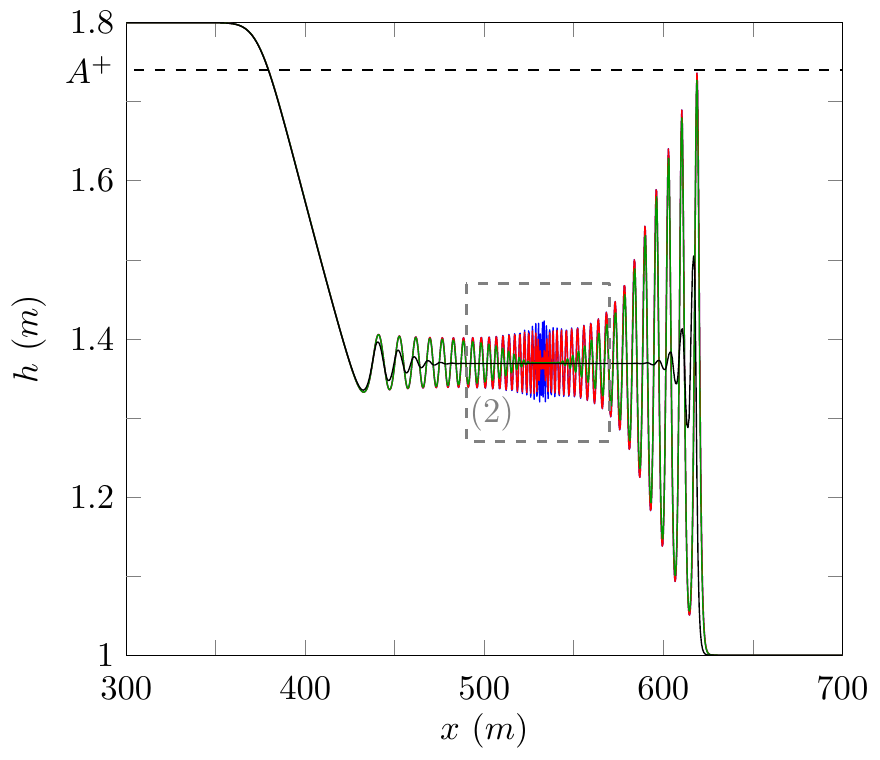}
\end{subfigure}
\begin{subfigure}{\textwidth}
	\includegraphics[width=0.5\textwidth]{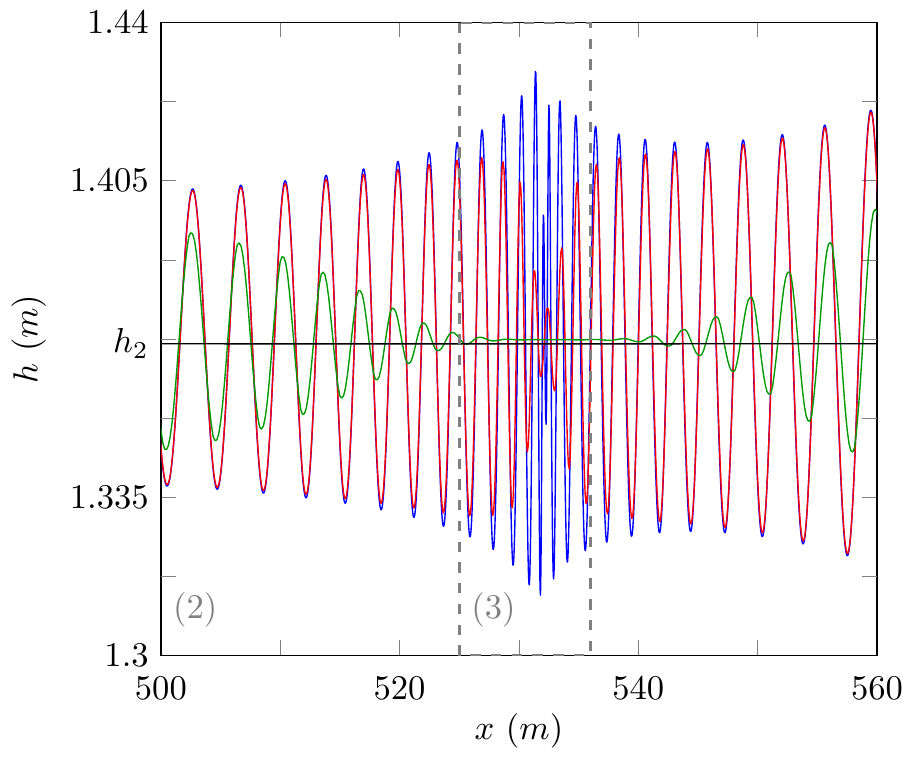}
	\includegraphics[width=0.5\textwidth]{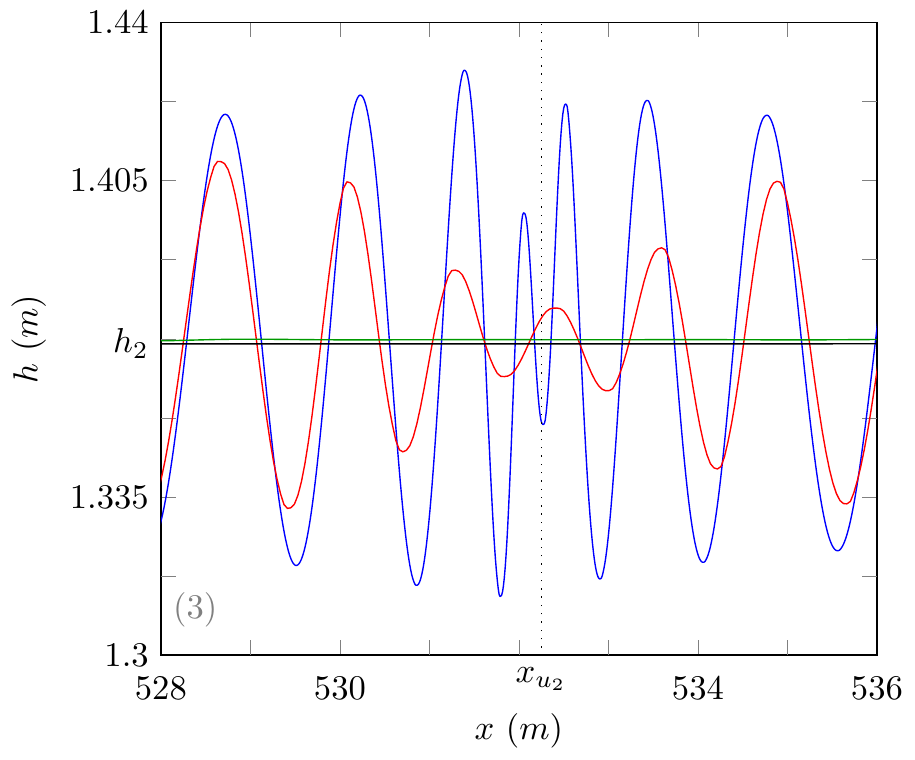}
\end{subfigure}
\caption{Numerical solutions of $\mathcal{V}_3$ at $t= 30s$ for the smooth dam-break problem with $\alpha = 0.1m$ for $\Delta x = 10/2^{10}m$ ({\color{blue} \solidrule}), $10/2^8m$ ({\color{red} \solidrule}), $10/2^6m$ ({\color{green!60!black} \solidrule}) and $10/2^{4}m$ ({\color{black} \solidrule}). The important quantities $A^+$ ({\color{black} \dashedrule}), $h_2$ ({\color{black} \dashedrule}) and $x_{u_2}$ ({\color{black} \dotrule{4mm}}) are also presented.}
\label{fig:o3a20dxlimcdexp}
\end{figure}

The $\mathcal{S}_4$ ``growth'' structure, which has hitherto not been commonly published in the literature features a growth in the oscillation amplitude around $x_{u_2}$. An example of the growth structure can be seen for $\mathcal{V}_3$'s numerical solutions in Figure \ref{fig:o3a20dxlimcdexp} to the smoothed dam-break problem with $\alpha = 0.1m$. This structure was observed in the numerical solutions of $\mathcal{V}_3$ for $\Delta x = 10/2^{11}m$ at $t=30s$ for $\alpha$ values as low as $0.001m$ and even for the dam-break problem.

Figure \ref{fig:o3a20dxlimcdexp} shows that this structure can only be observed for $\Delta x = 10 / 2^{10}m$, with poor convergence of the numerical results around $x_{u_2}$. Again our $L_1$ measures in Table \ref{tab:L1C1} omit the interval [$520m$,$540m$] in the numerical solutions. This demonstrates that although we have visual convergence away from $x_{u_2}$ our numerical solutions have not converged to round-off error as $\Delta x \rightarrow 0$. The $C_1$ measures in Table \ref{tab:L1C1} are still decreasing and have only attained round-off error for $h$, although for $uh$ and $\mathcal{H}$ the errors in conservation are small. These measures continue the trend in Table \ref{tab:L1C1} where smaller $\alpha$'s and thus steeper initial conditions lead to larger $L_1$ and $C_1$ measures because steeper problems are more difficult to solve accurately.

\begin{figure}
	\centering
	\begin{subfigure}{\textwidth}
		\includegraphics[width=0.49\textwidth]{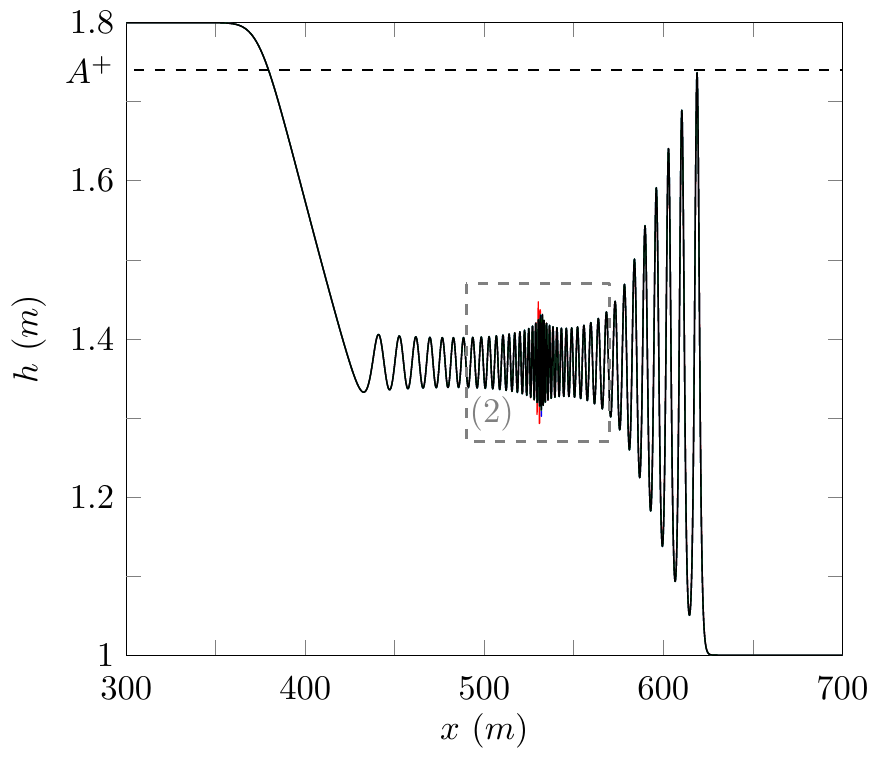}
		\includegraphics[width=0.5\textwidth]{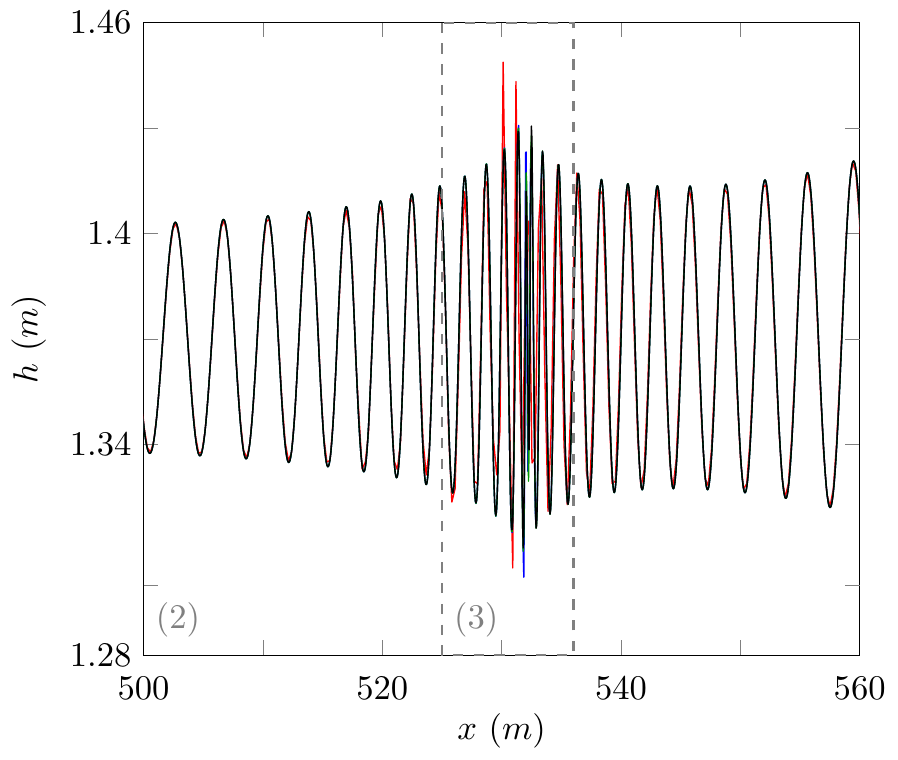}
	\end{subfigure}
	\begin{subfigure}{\textwidth}
		\includegraphics[width=0.5\textwidth]{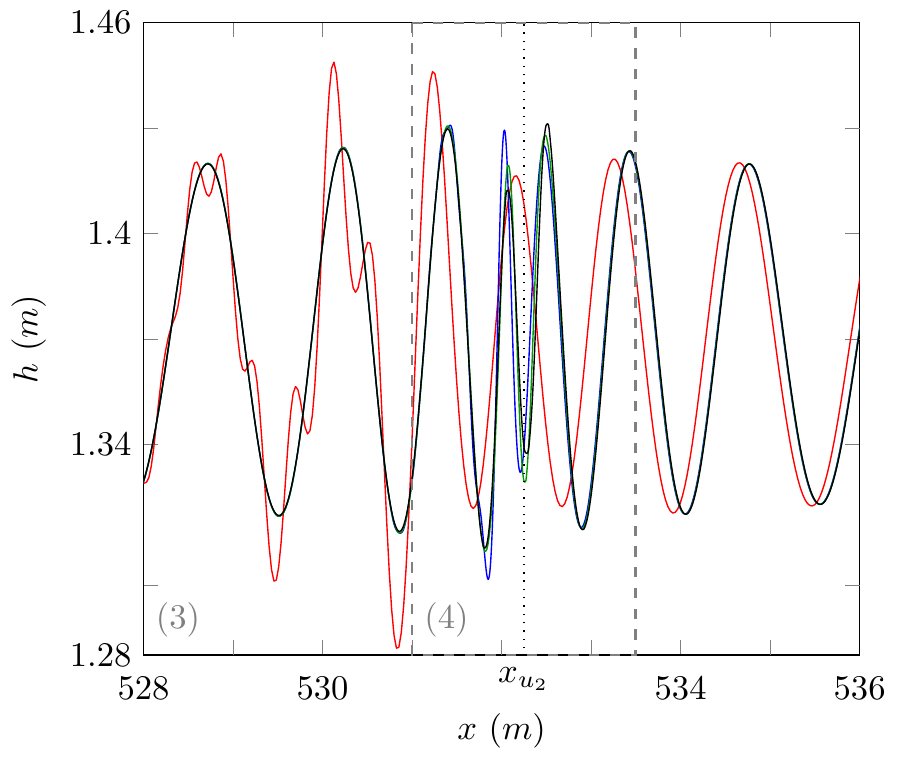}
		\includegraphics[width=0.5\textwidth]{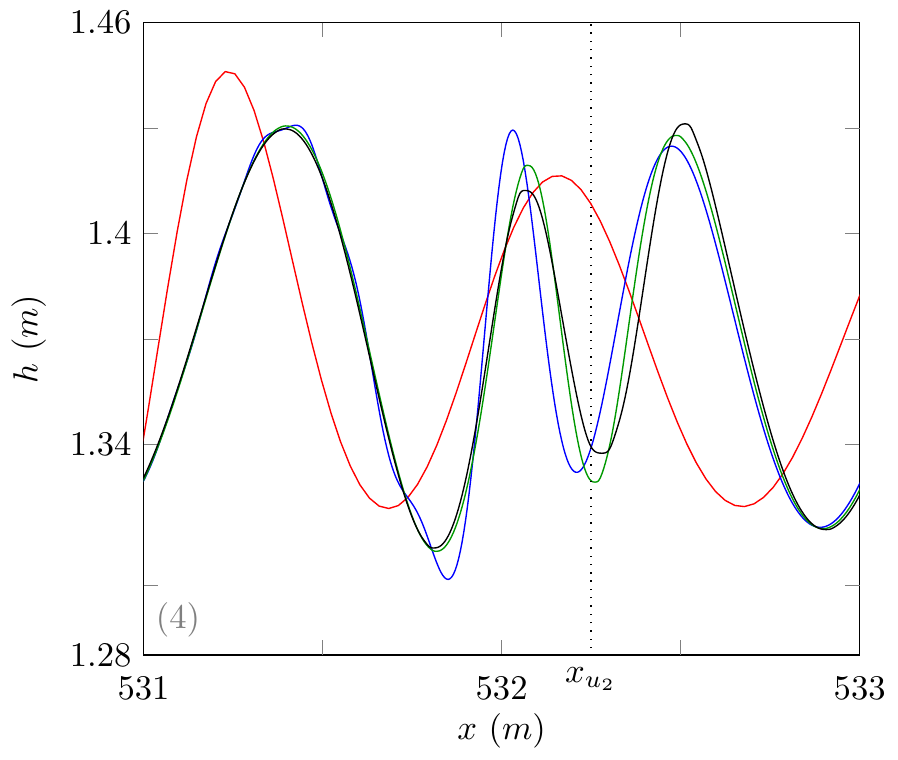}
	\end{subfigure}
	\caption{Numerical solutions of $\mathcal{D}$ ({\color{blue} \solidrule}), $\mathcal{E}$ ({\color{red} \solidrule}), $\mathcal{V}_3$ ({\color{green!60!black} \solidrule}) and $\mathcal{V}_2$ ({\color{black} \solidrule}) at $t=30s$ with $\Delta x = 10/2^{11}m$ for the smoothed dam-break problem with $\alpha = 0.1m$. The important quantities $A^+$ ({\color{black} \dashedrule}) and $x_{u_2}$ ({\color{black} \dotrule{4mm}}) are also presented.}
	\label{fig:o3a12allmodels}
\end{figure}

Figure \ref{fig:o3a12allmodels} demonstrates that our numerical solutions for $\Delta x = 10 /2^{11}m$ with the best methods $\mathcal{D}$, $\mathcal{V}_3$ and $\mathcal{V}_2$ disagree for only a few oscillations around $x_{u_2}$. Since both $\mathcal{D}$ and $\mathcal{E}$ are second-order finite difference methods their errors are dispersive. These dispersive errors cause the numerical solutions to overestimate the oscillation amplitude of the true solution, particularly around $x_{u_2}$. Because the dispersive errors of $\mathcal{E}$ are larger than $\mathcal{D}$ more oscillations are observed for the numerical solutions produced by $\mathcal{E}$. The $\mathcal{V}_3$ method was shown to be diffusive by \citet{Zoppou-etal-2017} and therefore its numerical solutions underestimate the oscillation amplitude in the true solution. Therefore, the true solution of the Serre equations should be between the dispersive method $\mathcal{D}$ and the diffusive method $\mathcal{V}_3$, and thus will possess the $\mathcal{S}_4$ structure.

\begin{figure}
	\centering
	\begin{subfigure}{\textwidth}
		\centering
		\includegraphics[width=0.5\textwidth]{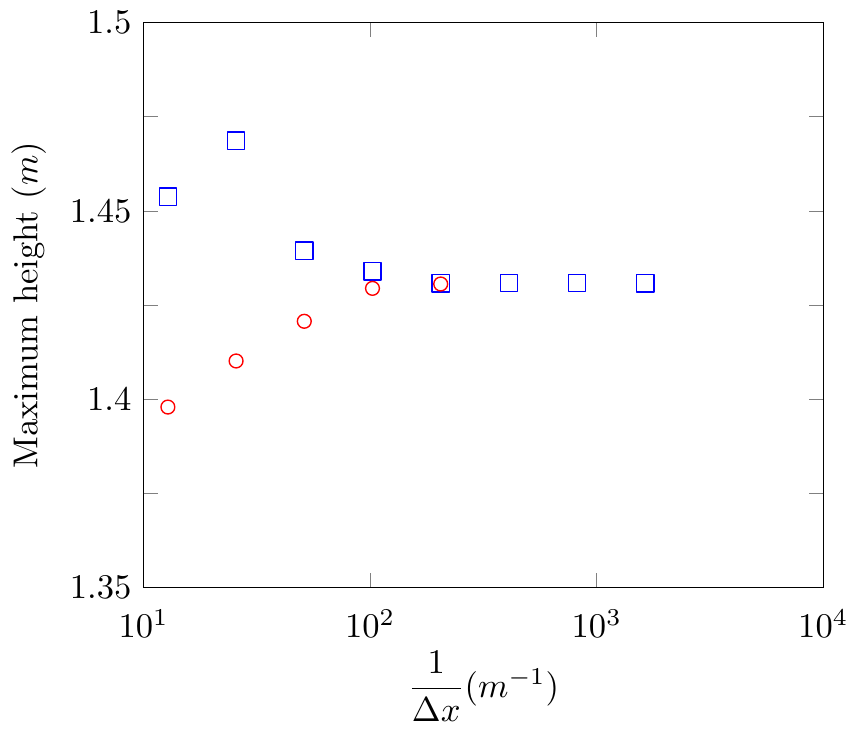}
	\end{subfigure}
	\caption{Maximum height of numerical solution of the smoothed dam-break problem with $\alpha = 0.4m$ at $t=30s$ inside the interval [$520m$,$540m$] using $\mathcal{D}$ ({\color{blue} $\square$}) and $\mathcal{V}_3$ ({\color{red} $\circ$}).}
	\label{fig:maxamp}
\end{figure}

The numerical solutions of $\mathcal{D}$ and $\mathcal{V}_3$ acting as upper and lower bounds respectively for the oscillation amplitude as $\Delta x$ is reduced is demonstrated in Figure \ref{fig:maxamp} using the maximum of $h$ in the interval [$520m$, $540m$]. From this figure it is clear that the amplitudes of the numerical solutions of $\mathcal{D}$ converge down  to the limit as the resolution is increased while the numerical solution amplitudes of $\mathcal{V}_3$ converge up to it. This shows that we have effectively bounded the true solution of the Serre equations. Unfortunately, $\mathcal{V}_3$ could not be run in reasonable computational times with lower $\Delta x$, but the numerical solutions of $\mathcal{D}$ show that doing so is unnecessary.

These results indicate that the solutions of the Serre equations to the smoothed dam-break problem with sufficiently small $\alpha$ values should exhibit a growth structure at $t=30s$, even though we have not precisely resolved all the oscillations in our numerical solutions. 

It was found that decreasing $\alpha$ did increase the amplitude of the oscillations around $x_{u_2}$. For $\mathcal{V}_3$ with $\Delta x= 10/2^{11}m$ and $\alpha = 0.001m$ the oscillations in $h$ were bounded by the interval [$1.28m$,$1.46m$]. Of particular note is that the number of oscillations are the same in Figures \ref{fig:o3a9allmodels} and \ref{fig:o3a12allmodels} for the best methods even though they have different structures.

\begin{figure}
	\centering
	\begin{subfigure}{\textwidth}
		\centering
		\includegraphics[width=0.5\textwidth]{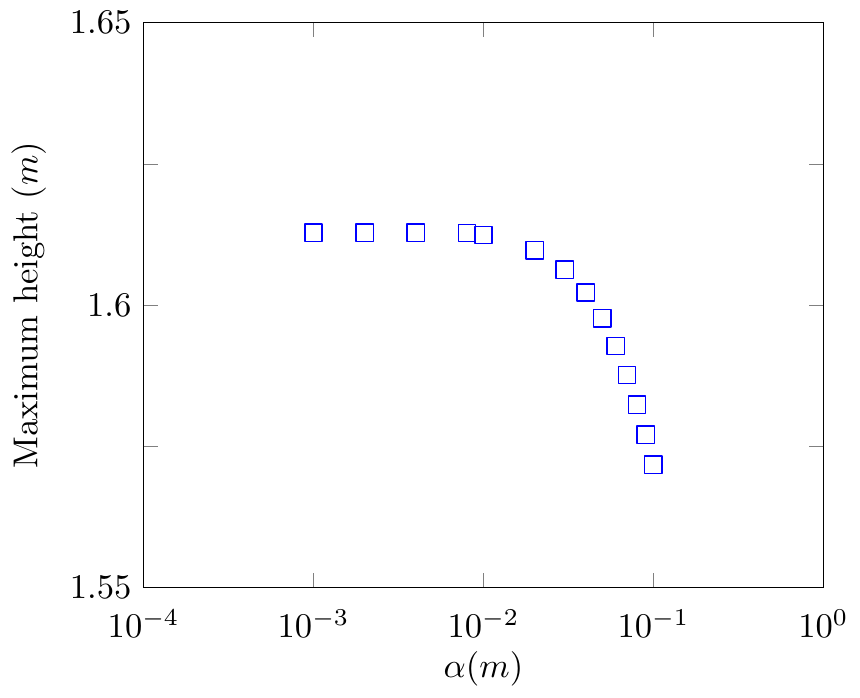}
	\end{subfigure}
	\caption{Maximum height of numerical solution around $x_{u_2}$ at $t=3s$ of various smoothed dam-break problem as $\alpha$ decreases, using $\mathcal{V}_3$ ({\color{blue} $\square$}) with $\Delta x = 10 / 2^{13}m$ .}
	\label{fig:maxampa}
\end{figure}

By changing the interval and desired time for the numerical solution, $\Delta x$ could be lowered further so that by $t=3s$ our numerical solutions have fully converged for $\alpha$ values as low as $0.001m$. This allows us to show that the height of the oscillations around $x_{u_2}$ for the solution of the Serre equation to the smoothed dam-break problem are bounded at $t=3s$ as $\alpha \rightarrow 0$. Figure \ref{fig:maxampa} demonstrates this for the numerical solutions of $\mathcal{V}_3$ with $\Delta x = 10/2^{13}m$.


\subsection{Shallow water wave equation comparison}
The analytical solutions of the shallow water wave equations have been used as a guide for the mean behaviour of the solution of the Serre equations for the dam-break problem in the literature \cite{Hank-etal-2010-2034,Mitsotakis-etal-2014}.

\begin{figure}
	\centering
	\begin{subfigure}{0.5\textwidth}
		\includegraphics[width=\textwidth]{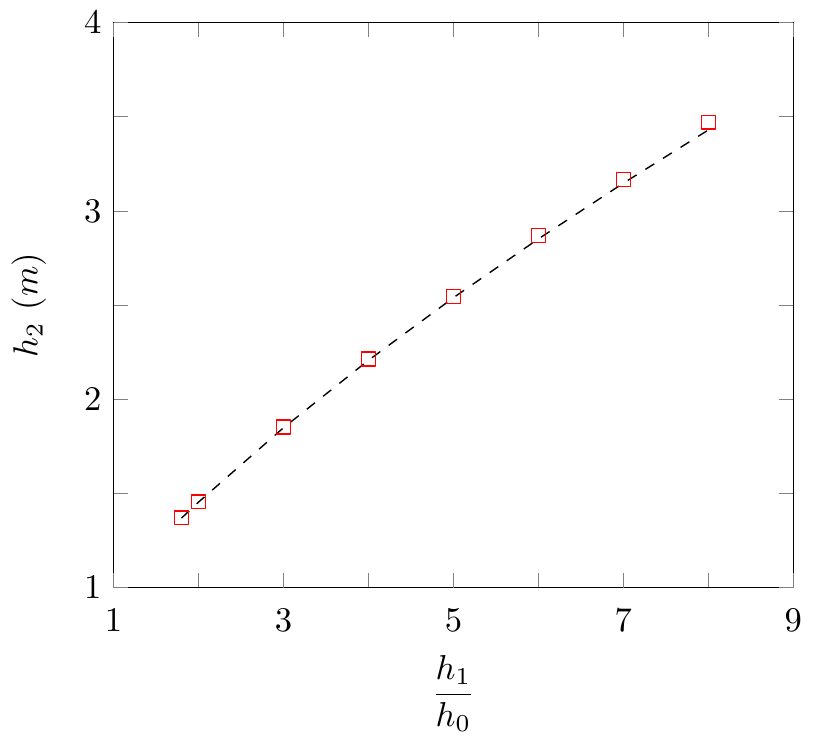}
		\subcaption*{\hspace{10 mm}h ({\color{red} $\square$}) and $h_2$({\color{black} \dashedrule})}
	\end{subfigure}%
	\begin{subfigure}{0.5\textwidth}
		\includegraphics[width=\textwidth]{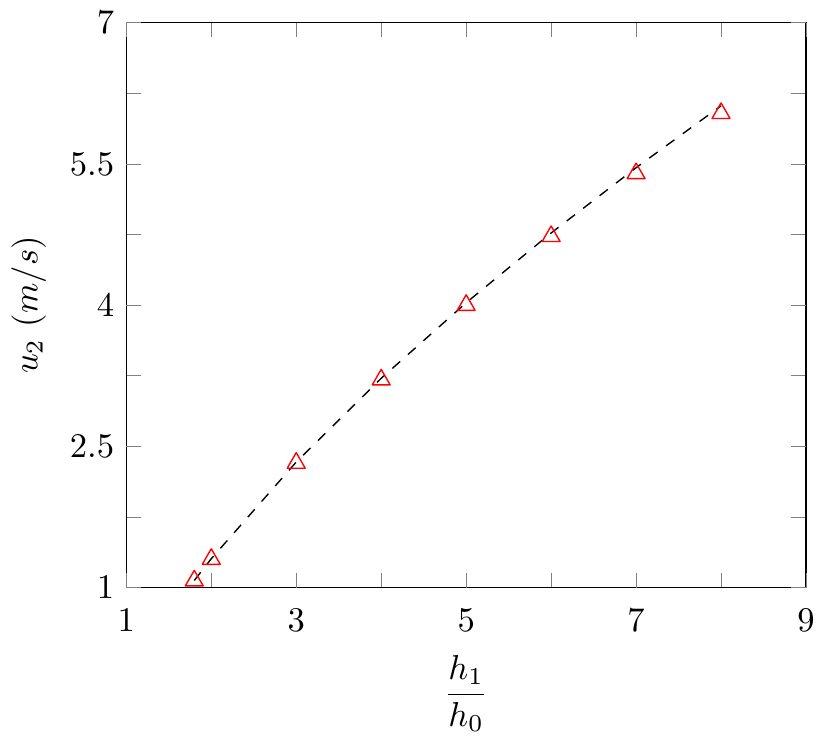}
		\subcaption*{\hspace{10 mm}u ({\color{red} $\triangle$}) and $u_2$({\color{black} \dashedrule})}
	\end{subfigure}
	\caption{Comparison between mean behaviour inside the bore of the Serre equations and the analytical solution of the shallow water wave equations for a range of different aspect ratios.}
	\label{fig:ratiotest}
\end{figure}

To assess the applicability of this the mean bore depth and mean fluid velocity in the interval [$x_{u_2}-50m$,$x_{u_2}+50m$] were calculated from our numerical solution to the smoothed dam-break problem with various height ratios. These means were compared to their approximations from the analytical solution of the dam-break problem for the shallow water wave equations $h_2$ and $u_2$. The results of this can be seen in Figure \ref{fig:ratiotest} for numerical solutions of $\mathcal{V}_3$ with $\Delta x = 10/2^{9}m$ to the smoothed dam-break problem at $t=100s$ with $\alpha = 0.1m$ where $h_0$ is fixed and $h_1$ is varied.

We use a final time of $t=100s$ as it allows the internal structure of the bore to develop more fully giving a more reliable mean estimate, as a consequence we resort to a coarser grid to keep the run-times reasonable. We find that decreasing $\Delta x$ does not significantly alter the mean of $h$ and $u$. We also find that increasing $\alpha$ also does not significantly alter the mean of $h$ and $u$. Therefore, the mean behaviour of the true solution of the Serre equations to the dam-break problem is captured by these numerical solutions, if it exists.

It can be seen that $h_2$ and $u_2$ are good approximations to the mean behaviour of the fluid inside the bore for a range of different aspect ratios. Although, as $h_1/h_0$ increases this approximation becomes worse, so that $h_2$ becomes an underestimate and consequently $u_2$ is an overestimate.

\begin{figure}
	\centering
	\begin{subfigure}{0.7\textwidth}
		\includegraphics[width=\textwidth]{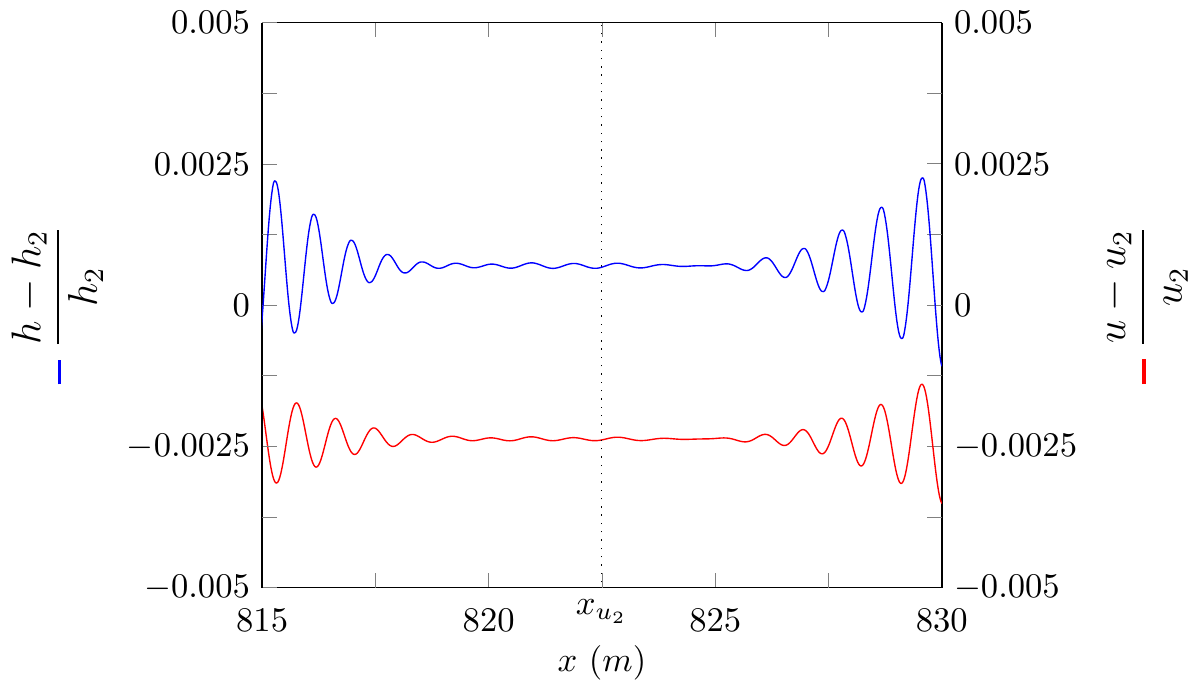}
	\end{subfigure}%
	\caption{The relative difference between $h$ and $u$ and their comparisons $h_2$ and $u_2$ plotted around $x_{u_2}$ ({\color{black} \dotrule{4mm}}) for $\mathcal{V}_3$'s solutions with $\Delta x = 10/2^{9}m$ for the smoothed dam-break problem with $\alpha = 0.1m$ at $t=300s$.}
	\label{fig:uhSWWcomp}
\end{figure}

We find that for $h_1/h_0 = 1.8$ the mean values of $h$ and $u$ inside the bore for the Serre equations are not equal to $h_2$ and $u_2$. This can be seen in Figure \ref{fig:uhSWWcomp} for the numerical solutions of $\mathcal{V}_3$ with $\Delta x = 10/2^{9}m$ to the smoothed dam-break problem with $\alpha = 0.1m$ at $t=300s$. It can be seen that $h_2$ is an underestimate of $h$ and $u_2$ is an overestimate of $u$ although the difference between these values and the mean behaviour of the Serre equations is small and only noticeable over long time periods.

\begin{figure}
	\centering
	\begin{subfigure}{0.5\textwidth}
		\includegraphics[width=\textwidth]{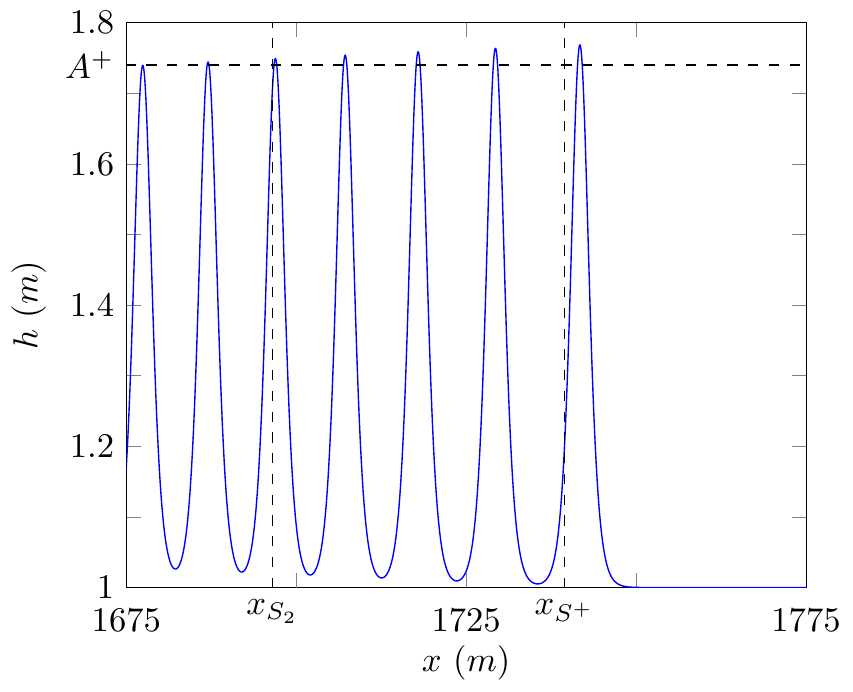}
	\end{subfigure}%
	\caption{Numerical solution of $\mathcal{V}_3$ with $\Delta x = 10/2^{9}m$ for the smoothed dam-break problem with $\alpha = 0.1m$ at $t=300s$ around the front of the undular bore. The important quantities $A^+$ ({\color{black} \dashedrule}), $x_{S_2}$ ({\color{black} \dashedrule}) and $x_{S^+}$ ({\color{black} \dashedrule}) are also presented.}
	\label{fig:frontcomp}
\end{figure}

The location of the leading wave of the Serre equations slowly diverges from the location of the front of a bore in the shallow water wave equations over long periods of time. This divergence causes the small difference evident in $\mathcal{V}_3$'s numerical solution to the smoothed dam-break problem with $\alpha =0.1m$ at $t=300s$ using $\Delta x = 10/2^{9}m$, which is shown in Figure \ref{fig:frontcomp}.

\begin{figure}
	\centering
	\begin{subfigure}{\textwidth}
		\centering
		\includegraphics[width=0.5\textwidth]{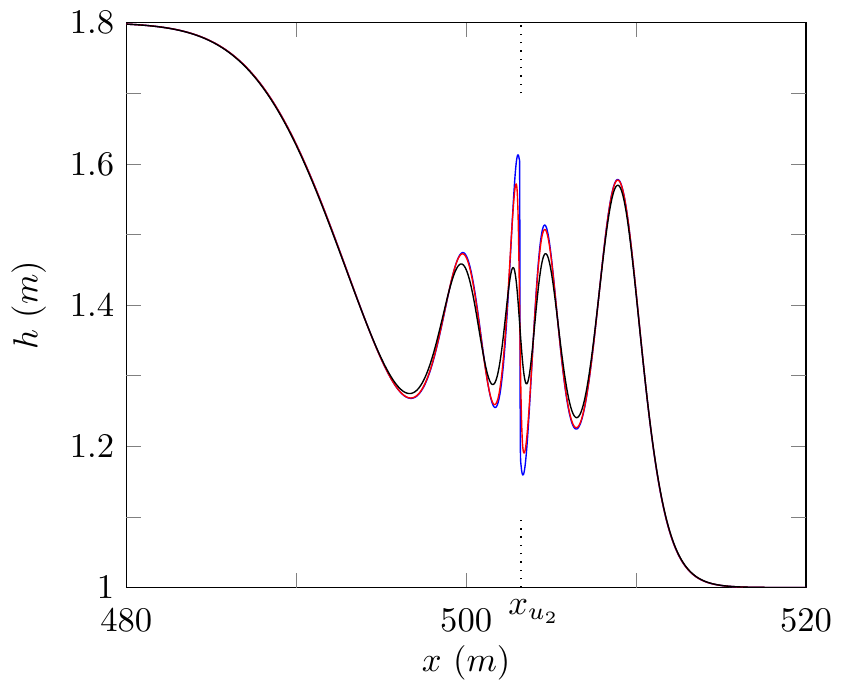}
	\end{subfigure}
	\caption{Numerical solution of $\mathcal{V}_3$ with $\Delta x = 10/2^{13}m$ for the smoothed dam-break problem with $\alpha = 0.001m$ ({\color{blue} \solidrule}), $0.1m$ ({\color{red} \solidrule}) and $0.4$ ({\color{black} \solidrule}) at $t=3s$. For comparison $x_{u_2}$ ({\color{black} \dotrule{4mm}}) is also plotted.}
	\label{fig:structsearlier}
\end{figure}

We note that the $\mathcal{S}_4$ structure present in the numerical solutions using this method and parameters at $t=30s$ in Figure \ref{fig:o3a20dxlimcdexp} has decayed away by $t=300s$ in Figure \ref{fig:uhSWWcomp}. This is a trend throughout our numerical solutions where oscillation amplitude decreases over time around $x_{u_2}$, changing the structure of the solution. This can be seen by obtaining full convergence of the numerical solutions to the smoothed dam-break problem at $t=3s$. The converged to numerical solutions for $\mathcal{V}_3$ are shown in Figure \ref{fig:structsearlier}. From this figure it can be seen that the oscillation amplitudes for the numerical solutions for the smoothed dam-break problems with $\alpha = 0.4m$ and $\alpha = 0.1m$ are much larger at $t=3s$ than they are at $t=30s$ in Figure \ref{fig:allstructs}. Since we have demonstrated that our numerical solutions are good approximations to the true solution of the Serre equations at $t=30s$ and $t=3s$, decreasing oscillation amplitude around $x_{u_2}$ over time is probably a property of the Serre equations. This implies that bounding the oscillation amplitudes at time $t=3s$ as was done above, bounds the oscillation amplitudes at all later times.
 
\subsubsection{Contact discontinuity}
\citet{El-etal-2006} noted the presence of a `degenerate contact discontinuity' which is the node in the $\mathcal{S}_3$ structure and travels at the mean fluid velocity in the bore.

We observe that as our numerical solutions evolve over time, oscillations appear to be released from the contact discontinuity and travel away from it in both directions, leading to decay of amplitudes around the contact discontinuity. Therefore, the contact discontinuity is an important feature and its behaviour determines the structure of the oscillations in the middle of the undular bore. 

The different speeds of the oscillations are determined by the phase velocity, which for the Serre equations linearised around the mean height $\bar{h}$ and mean velocity $\bar{u}$ in regions III and IV of the solution to the dam-break problem is
\begin{linenomath*}
	\begin{gather}
	\upsilon_p = \bar{u} \pm \sqrt{g\bar{h}} \sqrt{\frac{3}{\bar{h}^2 k^2 + 3}}
	\end{gather}
\end{linenomath*}
with wave number $k$. It can be seen that as $k \rightarrow \infty$ then $\upsilon_p \rightarrow \bar{u}$ and as $k \rightarrow 0$ then $\upsilon_p \rightarrow \bar{u} \pm \sqrt{g\bar{h}}$. Since the contact discontinuity travels at the mean velocity inside the bore, it corresponds to very high wave number oscillations. The oscillations on the left travel slower than the contact discontinuity and are therefore lower wave number oscillations associated with the phase velocity $ \bar{u} - \sqrt{g\bar{h}} \sqrt{3/ (\bar{h}^2 k^2 + 3)}$. The oscillations on the right travel travel quicker than the contact discontinuity and are therefore lower wave number oscillations associated with the phase velocity $ \bar{u} + \sqrt{g\bar{h}} \sqrt{3/ (\bar{h}^2 k^2 + 3)}$.

\begin{figure}
	\centering
	\begin{subfigure}{0.55\textwidth}\centering
		\includegraphics[width=\textwidth]{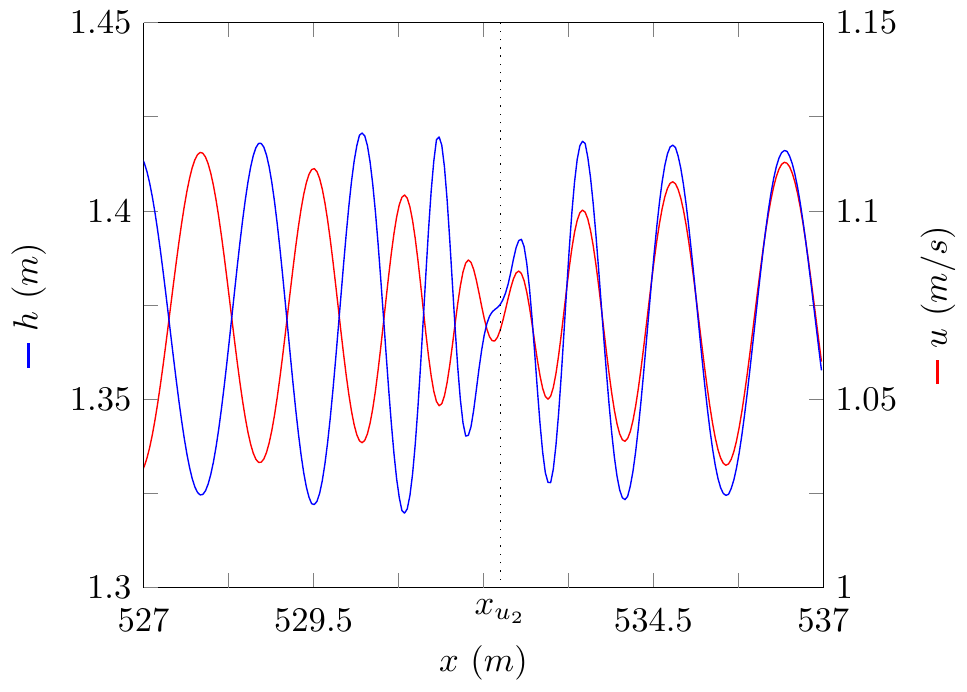}
	\end{subfigure}%
	\caption{Numerical solution of $\mathcal{V}_3$'s with $\Delta x = 10/2^{9}m$ for the smoothed dam-break problem with $\alpha = 0.1m$ at $t=30s$ around the contact discontinuity close to $x_{u_2}$ ({\color{black} \dotrule{4mm}}).}
	\label{fig:uhcomp}
\end{figure}

These different phase velocities have two different behaviours for $h$ and $u$. When the phase velocity is $ \bar{u} + \sqrt{g\bar{h}} \sqrt{3/ (\bar{h}^2 k^2 + 3)}$ we have oscillations where $h$ and $u$ are in-phase, while when the phase velocity is $ \bar{u} - \sqrt{g\bar{h}} \sqrt{3/ (\bar{h}^2 k^2 + 3)}$ we have oscillations where $h$ and $u$ are out-of-phase. This can be seen in Figure \ref{fig:uhcomp} for the numerical solutions of $\mathcal{V}_3$ with $\Delta x = 10/2^{9}m$ for the smoothed dam-break problem with $\alpha = 0.1m$ at $t=30s$.

\subsection{Whitham Modulation Comparsion}
\citet{El-etal-2006} demonstrated that their Whitham modulation results approximated the numerical solutions of the smoothed dam-break problem well for a range of aspect ratios. We observed that the Whitham modulation results are an underestimate compared to our numerical solutions.

\begin{figure}
	\centering
	\begin{subfigure}{0.5\textwidth}\centering
		\includegraphics[width=\textwidth]{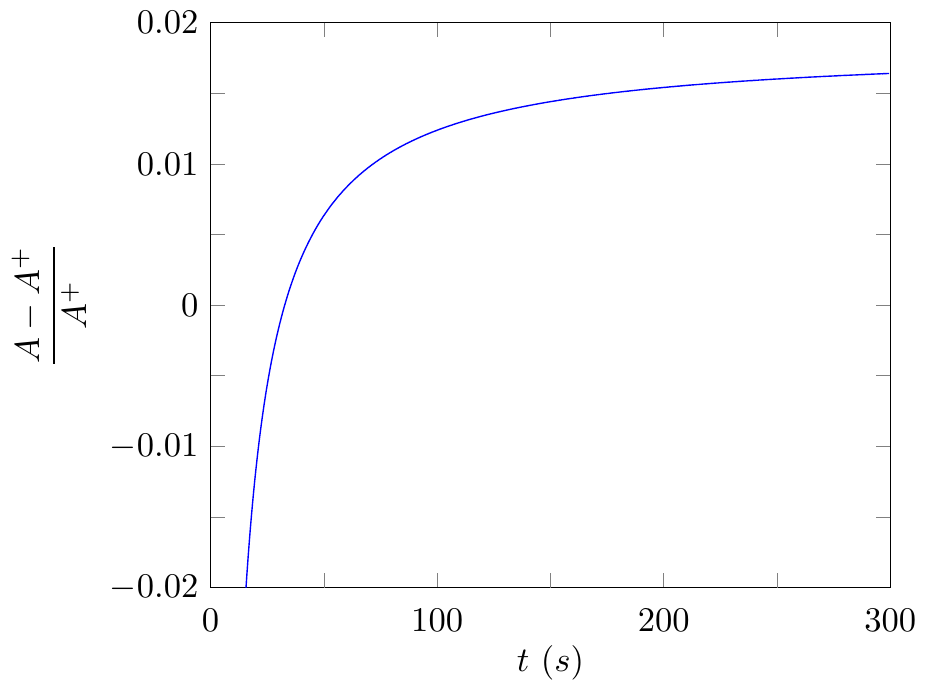}
	\end{subfigure}%
	\caption{Relative difference between Whitham modulation result $A^+$ and the leading wave amplitude $A$ from our numerical solutions of $\mathcal{V}_3$ with $\Delta x = 10/2^{9}m$ for the smoothed dam-break problem with $\alpha = 0.1m$ over time.}
	\label{fig:amplitudefront}
\end{figure}

This can be seen in Figure \ref{fig:amplitudefront} as the relative difference between $A^+$ from \citet{El-etal-2006} and the leading wave amplitude of our numerical solution $A$ does not converge to $0$ over time. Since we find that the numerical solutions for the smoothed dam-break problem with $\alpha = 0.1m$ have converged for the front of the undular bore by $\Delta x = 10/2^{8}m$ as in Figure \ref{fig:o3a20dxlimcdexp}, our numerical solutions for $A$ are considered reliable. We also note that unlike the oscillations around $x_{u_2}$ the leading wave amplitude increases over time.

The Whitham modulation results for the location of the leading wave $x_{S^+}$ is a better approximation than that given by the shallow water wave equations $x_{S_2}$, as can be seen in Figure \ref{fig:frontcomp}.
\section{Conclusions}
\label{section:Conclusions}
Utilising two finite difference methods of second-order and three finite difference finite volume methods of various orders to solve the nonlinear weakly dispersive Serre equations an investigation into the smoothed dam-break problem with varying steepness was performed. Four different structures of the numerical solutions were observed and demonstrated to be valid, the general trend of these structures is that an increase in steepness increases the size and number of oscillations in the solution. This study explains the different structures exhibited by the numerical results in the literature for the smoothed dam-break problem for the Serre equations and uncovers a new result. These results demonstrate that other methods in the literature could replicate our results if their simulations are extended. Furthermore, these results suggest that this new result and its associated structure is to be expected for the solution of the Serre equation to the dam-break problem at least for short enough time spans, if it exists.

We find that the analytical solution of the shallow water wave equations for the dam-break problem provides a reasonable approximation to the mean height and velocity inside the bore formed by the smoothed dam-break problem for the Serre equations. Finally, we observe that the Whitham modulations results for the leading wave of an undular bore provide a more accurate approximation to the location and depth of the front of an undular bore than the shallow water wave equations.


\section*{References}
\bibliography{DamBreak}
\bibliographystyle{elsarticle-num-names}
\newpage
\appendix{}
\label{sec:appendix}
\section{}
The methods $\mathcal{E}$ and $\mathcal{D}$ use the centred second-order finite difference approximation to the momentum equation \eqref{eq:Serre_momentum}, denoted as $\mathcal{D}_u$. For the mass equation \eqref{eq:Serre_continuity} $\mathcal{E}$ uses the two step Lax-Wendroff method, denoted as $\mathcal{E}_h$ while $\mathcal{D}$ uses a centred second-order finite difference approximation, denoted as $\mathcal{D}_h$.
\subsection{$\mathcal{D}_u$ for the Momentum Equation}
First \eqref{eq:Serre_momentum} is expanded to get
\begin{linenomath*}
		\begin{gather*}
		h\dfrac{\partial u}{\partial t} - h^2\frac{\partial^2 u}{\partial x \partial t} - \frac{h^3}{3}\frac{\partial^3 u}{\partial x^2 \partial t}  = -X 
		\label{eq:expandedu}
		\end{gather*}
		where $X$ contains only spatial derivatives and is
		\begin{gather*}
		X = uh\frac{\partial u}{\partial x} + gh\frac{\partial h}{\partial x} + h^2\frac{\partial u}{\partial x}\frac{\partial u}{\partial x} + \frac{h^3}{3}\frac{\partial u}{\partial x}\frac{\partial^2 u}{\partial x^2} - h^2u\frac{\partial^2 u}{\partial x^2}- \frac{h^3}{3}u\frac{\partial^3 u}{\partial x^3} .
		\end{gather*}
\end{linenomath*}
All derivatives are approximated by second-order centred finite difference approximations on a uniform grid in space and time, which after rearranging into an update formula becomes
\begin{linenomath*}
	\begin{gather}
	h^{n}_iu^{n+1}_i - \left(h^{n}_i\right)^2 \left(\frac{u^{n+1}_{i+1} -u^{n+1}_{i-1} }{2 \Delta x}\right) - \frac{\left(h^{n}_i\right)^3}{3}\left(\frac{u^{n+1}_{i+1} - 2u^{n+1}_{i} + u^{n+1}_{i-1} }{\Delta x^2}\right) = - Y^n_i
	\label{eq:expandedutdisc3}
	\end{gather}
\end{linenomath*}
where
\begin{linenomath*}
	\begin{gather*}
	Y_i^n = 2\Delta tX_i^{n} - h_i^{n}u_i^{n-1} + \left(h_i^{n}\right)^2\left(\frac{u^{n-1}_{i+1} -u^{n-1}_{i-1} }{2 \Delta x}\right) + \frac{\left(h_i^{n}\right)^3}{3}\left(\frac{u^{n-1}_{i+1} - 2u^{n-1}_{i} + u^{n-1}_{i-1} }{\Delta x^2}\right)
	\label{eq:expandfactor Xp}
	\end{gather*}
\end{linenomath*}
and
\begin{linenomath*}
\begin{multline*}
X_i^n = u_i^nh_i^n\frac{u^{n}_{i+1} -u^{n}_{i-1} }{2 \Delta x} + gh^n_i\frac{h^{n}_{i+1} -h^{n}_{i-1} }{2 \Delta x} + \left(h^n_i\right)^2\left(\frac{u^{n}_{i+1} -u^{n}_{i-1} }{2 \Delta x} \right)^2  \\ + \frac{\left(h^n_i\right)^3}{3}\frac{u^{n}_{i+1} -u^{n}_{i-1} }{2 \Delta x}\frac{u^{n}_{i+1} - 2u^{n}_{i} + u^{n}_{i-1} }{\Delta x^2} - \left(h^n_i\right)^2u_i^n\frac{u^{n}_{i+1} - 2u^{n}_{i} + u^{n}_{i-1} }{\Delta x^2} \\- \frac{\left(h^n_i\right)^3}{3}u^n_i \frac{u^{n}_{i+2} - 2u^{n}_{i+1} + 2u^{n}_{i-1} - u^{n}_{i-2} }{2\Delta x^3}.
\end{multline*}
\end{linenomath*}
 Equation \eqref{eq:expandedutdisc3} can be rearranged into an explicit update scheme $\mathcal{D}_u$ for $u$ given its current and previous values, so that
\begin{linenomath*}
	\begin{gather}
	\left[\begin{array}{c}
	u^{n+1}_0 \\
	\vdots \\
	u^{n+1}_m \end{array}\right]
	= A^{-1} \left[\begin{array}{c}
	-Y^n_0 \\
	\vdots \\
	-Y^n_m \end{array}\right] =: \mathcal{D}_u\left(\boldsymbol{u}^n,\boldsymbol{h}^n, \boldsymbol{u}^{n-1}, \Delta x, \Delta t \right)
	\label{eq:FDcentforu}
	\end{gather}
\end{linenomath*}
where $A$ is a tri-diagonal matrix.

\subsection{Numerical Methods for the Mass Equation}
\label{section:}
The two step Lax-Wendroff update $\mathcal{E}_h$ for $h$ is
\begin{linenomath*}
	\begin{gather*}
	h^{n + 1/2}_{i+ 1/2} = \frac{1}{2}\left(h^{n}_{i+1} + h^{n}_i\right) - \frac{\Delta t}{2\Delta x}\left(u^n_{i+1}h^n_{i+1} - h^n_{i}u^n_{i}\right),
	\end{gather*}
	\begin{gather*}
	h^{n + 1/2}_{i- 1/2} = \frac{1}{2}\left(h^{n}_{i} + h^{n}_{i-1}\right) - \frac{\Delta t}{2\Delta x}\left(u^n_{i}h^n_{i} - h^n_{i-1}u^n_{i-1}\right)
	\end{gather*}
	and
	\begin{gather*}
	h^{n+1}_i = h^{n}_i - \frac{\Delta t}{\Delta x}\left(u^{n + 1/2}_{i+ 1/2}h^{n + 1/2}_{i+ 1/2} - u^{n + 1/2}_{i- 1/2}h^{n + 1/2}_{i- 1/2}\right).
	\label{eq:LW4h}
	\end{gather*}
\end{linenomath*}
The quantities $u^{n + 1/2}_{i \pm 1/2}$ are calculated using $u^{n+1}$ obtained by applying $\mathcal{D}_u$ \eqref{eq:FDcentforu} to $u^n$ then linearly interpolating in space and time to give
\begin{linenomath*}
	\begin{gather*}
	u^{n + 1/2}_{i+ 1/2} = \frac{u^{n+1}_{i+1} + u^{n}_{i+1} + u^{n+1}_{i} + u^{n}_{i} }{4}
	\end{gather*}
	and
	\begin{gather*}
	u^{n + 1/2}_{i- 1/2} = \frac{u^{n+1}_{i} + u^{n}_{i} + u^{n+1}_{i-1}+ u^{n}_{i-1} }{4}.
	\end{gather*}
\end{linenomath*}
Thus we have the following update scheme $\mathcal{E}_h$ for \eqref{eq:Serre_continuity}
\begin{linenomath*}
	\begin{gather}
	\boldsymbol{h}^{n+1} = \mathcal{E}_h\left(\boldsymbol{u}^n,\boldsymbol{h}^n,\boldsymbol{u}^{n+1}, \Delta x, \Delta t \right). 
	\label{eq:LWupdateh}
	\end{gather}
\end{linenomath*}

The second order centered finite difference approximation to the conservation of mass equation \eqref{eq:Serre_continuity} is
\begin{linenomath*}
	\begin{gather*}
	h^{n+1}_i = h^{n-1}_i - \Delta t \left(u^{n}_{i}\frac{h^{n}_{i+1} - h^{n}_{i-1}}{\Delta x} + h^{n}_{i}\frac{u^{n}_{i+1} - u^{n}_{i-1}}{\Delta x}\right).
	\end{gather*}
\end{linenomath*}
Thus we have an update scheme $\mathcal{D}_h$ for all $i$
\begin{linenomath*}
	\begin{gather}
	\label{eq:secondFDappformass}
	\boldsymbol{h}^{n+1} = \mathcal{D}_h\left(\boldsymbol{u}^n,\boldsymbol{h}^n,\boldsymbol{h}^{n-1} ,\Delta x, \Delta t \right).
	\end{gather}
\end{linenomath*}

\subsection{Complete Method}
The method $\mathcal{E}$ is the combination of \eqref{eq:LWupdateh} for \eqref{eq:Serre_continuity} and \eqref{eq:FDcentforu} for \eqref{eq:Serre_momentum} in the following way
\begin{linenomath*}
	\begin{gather}
	\left.
	\begin{array}{l l}
	\boldsymbol{u}^{n+1}&=\mathcal{D}_u\left(\boldsymbol{u}^n,\boldsymbol{h}^n, \boldsymbol{u}^{n-1}, \Delta x, \Delta t \right) \\
	\boldsymbol{h}^{n+1}&=\mathcal{E}_h\left(\boldsymbol{u}^n,\boldsymbol{h}^n,\boldsymbol{u}^{n+1}, \Delta x, \Delta t \right)
	\end{array} \right\rbrace \mathcal{E}\left(\boldsymbol{u}^n,\boldsymbol{h}^n, \boldsymbol{u}^{n-1},\boldsymbol{h}^{n-1}, \Delta x, \Delta t \right).	 
	\end{gather}
\end{linenomath*}

The method $\mathcal{D}$ is the combination of \eqref{eq:secondFDappformass} for \eqref{eq:Serre_continuity} and \eqref{eq:FDcentforu} for \eqref{eq:Serre_momentum} in the following way
\begin{linenomath*}
	\begin{gather}
	\left.
	\begin{array}{l l}
	\boldsymbol{h}^{n+1}&=\mathcal{D}_h\left(\boldsymbol{u}^n,\boldsymbol{h}^n,\boldsymbol{h}^{n-1} \Delta x, \Delta t \right) \\
	\boldsymbol{u}^{n+1}&=\mathcal{D}_u\left(\boldsymbol{u}^n,\boldsymbol{h}^n, \boldsymbol{u}^{n-1}, \Delta x, \Delta t \right)
	\end{array} \right\rbrace \mathcal{D}\left(\boldsymbol{u}^n,\boldsymbol{h}^n, \boldsymbol{u}^{n-1},\boldsymbol{h}^{n-1}, \Delta x, \Delta t \right).
	\label{eq:Gnumdef}
	\end{gather}
\end{linenomath*}

\end{document}